\documentclass[11pt]{article}
\usepackage{mathptmx}       
\usepackage{fullpage} 
\usepackage{booktabs}

\usepackage{makeidx}         
\usepackage{graphicx}

\usepackage{multicol}        
\usepackage[bottom]{footmisc}

\usepackage{tcolorbox}

\usepackage{newtxtext}      
\usepackage[varvw]{newtxmath}       

\usepackage{graphicx} 
\usepackage[margin=1in]{geometry}
\usepackage{url}
\usepackage{tikz, comment}
\newcommand{\defeq}{\mathrel{\mathop:}=}

\newcommand{\iu}{\mathrm{i}\mkern1mu}

\usepackage{xcolor}

\usetikzlibrary{positioning, arrows.meta, calc, fit, backgrounds}
\newcommand{\bfx}{\mathbf{x}}
\newcommand{\bfy}{\mathbf{y}}
\newcommand{\bfz}{\mathbf{z}}
\newcommand{\bfd}{\mathbf{d}}
\newcommand{\bfs}{\mathbf{s}}
\newcommand{\bfv}{\mathbf{v}}

\title{Neural network enhancement of the Born approximation}
\author{J. Ma\thanks{Department of Mathematical Sciences, University of Delaware, Newark DE 19716, USA, {johnma@udel.edu}} \and A. Desai\thanks{Department of Mathematics, University of Chicago,  Chicago IL 60637, USA, {ansh@uchicago.edu}} \and T. L\"ahivaara\thanks{Department of Technical Physics, University of Eastern Finland, Kuopio 70211, Finland, {timo.lahivaara@uef.fi}} \and P. Monk\thanks{Department of Mathematical Sciences, University of Delaware, Newark DE 19716, USA, {monk@udel.edu}, ORCID: 0000-0002-6539-5897.}}

\begin{document}
\maketitle
\abstract{
We examine the use of neural networks in conjunction with the Born approximation for solving the inverse scattering problem for a penetrable medium in the frequency domain.  The goal is a fast method for determining the refractive index of a scatterer from multistatic far-field measurements at one frequency. We present a new neural network model that seeks to improve on the results in
Ref.~\cite{DMLM26} by sharpening the reconstructed images. In addition, we present two other networks motivated
by the corrected Born method that avoid the need to
invert the regularized Born operator.  These latter algorithms will be important for future three dimensional applications when the inverse Born approximation is difficult to compute, or when near-field measurements are used. Computational results illustrate that the new corrections to the Born approximation can provide an accurate and flexible approach to the inverse scattering problem. 
}

\section{Introduction}

Inverse scattering seeks to recover the material properties or geometry of an inaccessible object from remote measurements of scattered-field data. Problems of this type arise in fields such as acoustics, radar, seismology, and biomedical imaging, and they have long been studied within the mathematical inverse problems community, see, e.g., \cite{colton+kress4,CCH,Kirsch}. In this work, we consider a model inverse scattering problem for time-harmonic acoustic waves, where the goal is to estimate the refractive index of a penetrable scatterer from  multistatic far-field data. Our main focus is on combining model-based inversion with modern neural-network architectures in a way that preserves useful physical structure while improving robustness and reconstruction quality beyond the classical weak-scattering regime.

In the frequency-domain formulation, the total field satisfies the Helmholtz equation with a Sommerfeld radiation condition imposed on the scattered field at infinity. The far field pattern of the scattered field -then gives the asymptotic amplitude of the outgoing scattered wave. - The corresponding inverse map from far field data to the refractive index is ill-posed and requires regularization \cite{colton+kress4,Devaney:12}. In practical imaging systems, measurement noise, modeling error, limited aperture, and sparse receiver coverage further degrade reconstruction quality. This is especially relevant in biomedical imaging applications such as ultrasound tomography and full-wave breast imaging, where large amounts of scattered-field data are used to recover quantitative information about tissue  \cite{taskin21,lucka22,crtm23}.

A classical approach to remove the nonlinearity of inverse scattering is the Born approximation. In the forward problem, the Born approximation arises by truncating the Neumann series representation of the total field at zeroth order, thereby replacing the full nonlinear scattering map by a linearized operator \cite{Hudson:81,Born_Wolf:59,Devaney:82,Moskow:19}. For 
sufficiently weak scattering, this leads to an efficient approximate inversion procedure after regularization. The regularized inverse Born approximation remains attractive because it provides a transparent physical interpretation, admits fast implementations, and often gives a reasonable first estimate of the scatterer. However, its validity deteriorates rapidly when the contrast is high, multiple scattering is significant, or the imaging setup departs from the idealized regime in which the approximation is justified. In such cases, reconstructions based solely on Born inversion typically  estimate the contrast inaccurately and lose structural detail; see \cite{Slaney}.

Recent years have seen a rapid expansion of data-driven methods for inverse scattering and related wave-imaging problems. Broadly, two general approaches can be distinguished. In the first, a neural network is trained to approximate the forward problem, so that the fast learned model can be used within a traditional optimization framework for the inverse problem; see \cite{Stanziola_2021,Chen:24,liu2024neumannseriesbasedneuraloperator}. In the second, end-to-end algorithms are developed that seek to solve the inverse problem directly by learning the map from scattering data to the unknown coefficient; see, for instance, \cite{lexing2019,gao21,Fan:22,ZHANG24,Bazow,NIO}. More particularly, operator-learning frameworks such as DeepONet, neural operators, and the Fourier Neural Operator have expanded the range of architectures available for PDE-based forward and inverse problems; see, for instance, \cite{Lu_2021,DBLP:journals/corr/abs-2108-08481,li2020fourier}.

Nevertheless, purely data-driven end-to-end approaches for inverse scattering often exhibit familiar weaknesses: they may require large training sets, be sensitive to distribution shift, and hallucinate structures when tested on noisy or out-of-distribution data. These issues are now well recognized more broadly in learned inverse problems \cite{siam-rev}. For inverse scattering, such behavior is especially problematic because small artifacts can be mistaken for physically meaningful inclusions.

These weaknesses have motivated hybrid strategies in which approximate physical models are not discarded, but instead used to provide additional physical information to the network, for example, as a model constraint, intermediate layer, or corrective operator. A related philosophy appears in model-guided or model-based learning for inverse problems and imaging  \cite{Arridge2019,hauptman18,hoskins25}. Closely related learning-assisted and physics-embedded approaches have also been proposed for nonlinear and full-wave inverse scattering, combining approximate inversion procedures, forward modeling, or iterative data-consistency mechanisms with trainable networks  \cite{9996181,9509361}.  Of particular relevance to this paper are algorithms that correct an approximate inversion model using a neural network; for example, in \cite{DeepNIS2019} 
neural networks are used to improve the efficiency of iterative methods for inverse scattering, while in \cite{WeiChen2019} a U-Net is used to correct back-propagation methods. By contrast, in  \cite{Jin2017} neural networks are used to correct a direct inversion of  CT data resulting in a general inversion scheme called FBPConvNet. Physics-informed neural networks and related PDE-constrained learning methods represent a distinct line of work in which the governing equations are incorporated directly into the training procedure and therefore do not typically depend on large paired datasets, although they can face optimization difficulties and involve per-instance training cost, see, e.g., \cite{RAISSI2019686,yuyao_2019,yuyao_2022,Saba,wu2022}.

Recent hybrid operator-learning work has also emphasized that high-contrast wave scattering is an especially challenging regime, and that decomposing smooth background propagation from scattering by jumps in the contrast can improve surrogate accuracy \cite{balaji2026hybridoperatorlearningwave}. For broader overviews, see, e.g.,  \cite{ChenReview2020,lin2025survey,maokun26}, which emphasize both the value of learned inversion and the continuing importance of physics-based priors and hybrid methods. Previous work has shown that errors due to approximate forward models can be reduced using statistical or learned model-error compensation approaches in scattering and acoustic inversion \cite{Kaipio:19,koponen}.

Our earlier work \cite{DMLM26} explored this general direction by introducing CNN-based correction strategies for Born inversion in the far-field setting. Those results indicated that embedding the Born approximation into a learning pipeline can improve reconstruction quality beyond the weak-scattering regime while preserving some interpretability from the underlying physics. In particular, they suggested that classical linearizations such as the Born approximation, when combined with modern neural components, may provide a useful physics-motivated structure that guides learning for inverse scattering.

Our main goal for this paper is to investigate three new networks
for the inverse scattering problems in Section~\ref{sec:New_NN}.  The motivation for the new networks 
is twofold: 1) we want to improve the resolution of the output of the 
inverse solver (without increasing hallucination unacceptably) and 2) prepare for the use of our methods for 3D near-field problem such as those arising in medical imaging~\cite{lucka22}.  In the latter case the Born approximation ceases to be
a band limited Fourier transform and so its inversion slows down.  

Besides testing the above mentioned networks for solving inverse scattering problems, we shall investigate  how our 
methods perform on impenetrable scatterers.  Since our networks are trained only on penetrable scatterers, this is a strong test of generalization. Secondly we will test how our networks perform on near-field data.  In this case we retrain all the networks on data measured
approximately one wavelength from the inspection domain where we a priori know the scatterers to lie.  In this case the Born approximation must also be modified.  Preliminary tests are performed on in-distribution scatterers.

The paper proceeds as follows: in Section~\ref{sec:forward-inverse} we 
present the forward and inverse problems under study and indicate how the Born approximation arises.
Then, in Section~\ref{sec:New_NN}, we also briefly discuss the neural network-corrected Born approximation networks in~\cite{DMLM26} and present details of the new neural networks under study.  Numerical results for the networks, all trained on the same synthetic data, are given in Section~\ref{sec:perf}.  This starts with a study of how the networks handle added noise on the 
data.  This gives information about the degree of hallucination, as well
as stability of the algorithms, and is discussed more in Section~\ref{sec:noise}.  We then consider two out-of-distribution cases: first the reconstruction of the out of distribution
scatterers from \cite{DMLM26}, and then reconstruction of impenetrable obstacles.  The networks in this paper are all trained on data from penetrable scatterers.  We are interested to determine how the networks
respond if the scattering mechanism is incorrectly modeled.  So we apply the networks to data from sound hard scatterers in Section~\ref{sec:obstacles}.  Next, in Section~\ref{sec:near}, we start the study of inversion of near-field scattered data.  This requires training new models on near-field data.  Finally, in Section~\ref{sec:conc} we draw some conclusions from our study.

\section{Problem setting}
\label{sec:forward-inverse}

In this section we summarize the  forward and inverse problems. We study  single frequency acoustic scattering governed by the Helmholtz equation in two space dimensions. We assume $\exp(-i\omega t)$ dependence on time.  Here $\omega$ is the angular frequency of the wave and $i=\sqrt{-1}$. Bold-face symbols represent vectors.
 We shall consider two types of scatterers: 1)  penetrable objects in which the wave field couples inside the scatterer, and 2) impenetrable obstacles with a Neumann  (or sound hard) boundary condition.  

\subsection{Forward problem}
We start by describing the standard penetrable scattering problem that was used in \cite{DMLM26}. Let $n(\bfx)$ denote the spatially varying refractive index of the medium, with $\bfx \in \mathbb{R}^2$. We assume $n \in L^{\infty}(\mathbb{R}^2)$ is piecewise smooth and that there exist constants $n_{\min}$ and $n_{\max}$ such that
\[
0 < n_{\min} \leq n(\bfx) \leq n_{\max} < \infty
\quad \text{for almost every } \bfx \in \mathbb{R}^2.
\]
The scatterer is assumed bounded, so that if $\|\bfx\| > R$ for some sufficiently large $R > 0$, then $n(\bfx) = 1$. We denote the support of the scatterer by the set $D$.

A known incident field impinges on the scatterer and results in a scattered field.  The incident field, denoted $u^{i}$, is assumed to satisfy the background Helmholtz equation $\Delta u^{i}+\kappa^2u^{i}=0$ in the neighborhood of the scatterer $D$ where $\kappa = \omega / c_0$ is the wave number in the background medium and $c_0$ the sound speed of waves outside $D$. 

For most of the paper we use plane incident waves given by
\begin{equation}\label{eq:ui_plane}
    u^{\mathrm{i}}(\bfx;\mathbf{d}) = \exp(i\kappa\, \bfd \cdot \bfx), 
    \qquad \|\bfd\| = 1,
\end{equation}
where the unit vector $\bfd$ specifies the direction of propagation of the wave.  This is consistent with our assumption of far-field measurements for  the inverse problem.

When we consider near-field scattering, we assume the incident field is due to a point source at a point $\bfz\not\in \overline{D}$ (where $\overline{D}$ is the closure of $D$) and is given by
\begin{equation}\label{eq:ui_point}
    u^{\mathrm{i}}(\bfx;\mathbf z) = \frac{i}{4} H^{(1)}_0(\kappa \Vert \bfx -\bfz\Vert),\end{equation}
where $H_0^{(1)}$ is the Hankel function of the first kind and zero order.

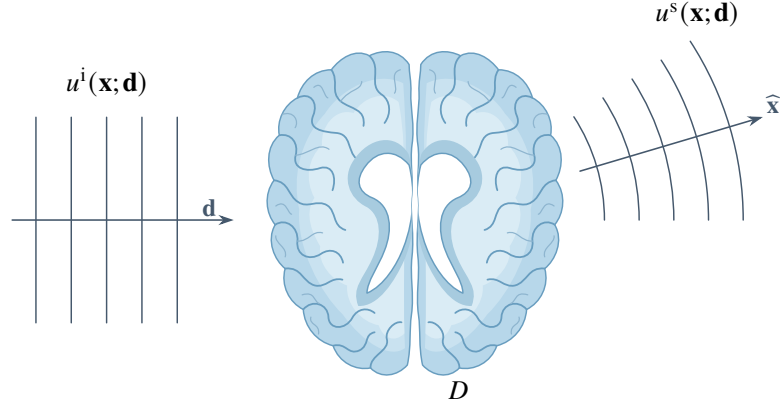
\begin{figure}
    \centering
\resizebox{0.7\textwidth}{!}{\begin{tikzpicture}[
  x=1cm,y=1cm,line cap=round,line join=round,
  >={Stealth[length=2.4mm,width=1.6mm]},
  field line/.style={line width=0.65pt},
  direction/.style={->,line width=0.7pt},
  guide/.style={draw=black!55,line width=0.45pt},
  sulcus/.style={draw=brainfold,line width=0.80pt},
  fine sulcus/.style={draw=brainfold!75,line width=0.55pt},
  every node/.style={inner sep=1.5pt}
]
\definecolor{incident}{RGB}{69,88,108}
\definecolor{scattered}{RGB}{69,88,108}
\definecolor{braintissue}{RGB}{183,215,234}
\definecolor{brainlight}{RGB}{218,236,245}
\definecolor{brainmid}{RGB}{201,226,240}
\definecolor{brainfold}{RGB}{104,157,191}
\definecolor{brainrim}{RGB}{145,187,212}
\path[use as bounding box] (0,0) rectangle (13.6,7.2);
\begin{scope}[shift={(6.55,3.65)},scale=0.91]
  \def\hemisphere{
    (-0.10,2.66)
    .. controls (-0.27,2.92) and (-0.56,2.92) .. (-0.78,2.81)
    .. controls (-1.02,2.94) and (-1.28,2.77) .. (-1.43,2.57)
    .. controls (-1.73,2.57) and (-1.94,2.35) .. (-2.01,2.09)
    .. controls (-2.24,1.97) and (-2.37,1.70) .. (-2.32,1.45)
    .. controls (-2.55,1.26) and (-2.54,0.91) .. (-2.41,0.68)
    .. controls (-2.59,0.44) and (-2.55,0.11) .. (-2.39,-0.08)
    .. controls (-2.55,-0.35) and (-2.40,-0.69) .. (-2.22,-0.83)
    .. controls (-2.30,-1.12) and (-2.15,-1.41) .. (-1.91,-1.53)
    .. controls (-1.89,-1.83) and (-1.68,-2.08) .. (-1.40,-2.15)
    .. controls (-1.30,-2.44) and (-1.01,-2.58) .. (-0.76,-2.53)
    .. controls (-0.58,-2.72) and (-0.30,-2.67) .. (-0.13,-2.48)
    .. controls (-0.06,-2.15) and (-0.10,-1.90) .. (-0.06,-1.62)
    .. controls (-0.01,-1.26) and (-0.09,-0.97) .. (-0.05,-0.62)
    .. controls (-0.01,-0.25) and (-0.08,0.05) .. (-0.06,0.37)
    .. controls (-0.02,0.69) and (-0.09,0.99) .. (-0.05,1.29)
    .. controls (0.00,1.62) and (-0.08,1.91) .. (-0.07,2.18)
    .. controls (-0.06,2.37) and (-0.10,2.52) .. (-0.10,2.66)
    -- cycle
  }
  \foreach \mirror in {1,-1} {
    \begin{scope}[xscale=\mirror]
      \path[fill=braintissue,draw=brainfold,line width=0.65pt] \hemisphere;
      \begin{scope}
        \clip \hemisphere;
        \path[fill=brainmid]
          (-0.20,2.34)
          .. controls (-0.67,2.58) and (-1.29,2.22) .. (-1.54,1.82)
          .. controls (-1.98,1.35) and (-2.10,0.72) .. (-1.98,0.18)
          .. controls (-2.03,-0.47) and (-1.61,-1.26) .. (-1.24,-1.61)
          .. controls (-0.87,-2.12) and (-0.42,-2.38) .. (-0.17,-2.19)
          .. controls (-0.22,-1.11) and (-0.18,0.01) .. (-0.17,0.81)
          .. controls (-0.13,1.41) and (-0.18,1.96) .. (-0.20,2.34) -- cycle;
        \path[fill=brainlight]
          (-0.28,2.02)
          .. controls (-0.53,2.23) and (-0.89,2.10) .. (-1.05,1.83)
          .. controls (-1.37,1.93) and (-1.62,1.58) .. (-1.59,1.28)
          .. controls (-1.83,1.06) and (-1.89,0.69) .. (-1.72,0.38)
          .. controls (-1.91,0.09) and (-1.74,-0.25) .. (-1.62,-0.45)
          .. controls (-1.67,-0.75) and (-1.42,-1.03) .. (-1.16,-1.12)
          .. controls (-1.19,-1.43) and (-0.92,-1.75) .. (-0.62,-1.80)
          .. controls (-0.55,-2.02) and (-0.35,-2.11) .. (-0.18,-1.94)
          .. controls (-0.20,-1.24) and (-0.17,-0.39) .. (-0.19,0.32)
          .. controls (-0.17,0.96) and (-0.19,1.59) .. (-0.28,2.02) -- cycle;
        \draw[sulcus]
          (-0.78,2.81) .. controls (-0.55,2.59) and (-0.82,2.31) .. (-0.56,2.14)
          .. controls (-0.42,2.04) and (-0.40,1.82) .. (-0.52,1.66);
        \draw[sulcus]
          (-1.43,2.57) .. controls (-1.29,2.37) and (-1.46,2.20) .. (-1.27,2.03)
          .. controls (-1.09,1.89) and (-0.84,1.91) .. (-0.88,1.61);
        \draw[fine sulcus]
          (-1.27,2.03) .. controls (-1.12,2.25) and (-0.92,2.21) .. (-0.91,2.37);
        \draw[sulcus]
          (-2.01,2.09) .. controls (-1.79,2.20) and (-1.65,2.01) .. (-1.74,1.79)
          .. controls (-1.84,1.55) and (-1.58,1.51) .. (-1.35,1.49)
          .. controls (-1.08,1.47) and (-1.16,1.18) .. (-1.00,1.05);
        \draw[sulcus]
          (-2.32,1.45) .. controls (-2.11,1.49) and (-1.99,1.31) .. (-2.04,1.11)
          .. controls (-2.10,0.88) and (-1.86,0.81) .. (-1.62,0.89)
          .. controls (-1.40,0.95) and (-1.34,0.72) .. (-1.41,0.51);
        \draw[fine sulcus]
          (-2.04,1.11) .. controls (-1.82,1.12) and (-1.65,1.34) .. (-1.47,1.21);
        \draw[sulcus]
          (-2.41,0.68) .. controls (-2.15,0.65) and (-2.17,0.41) .. (-2.00,0.31)
          .. controls (-1.84,0.19) and (-1.67,0.36) .. (-1.50,0.25)
          .. controls (-1.34,0.15) and (-1.42,-0.10) .. (-1.19,-0.17);
        \draw[sulcus]
          (-2.39,-0.08) .. controls (-2.21,0.04) and (-1.99,-0.12) .. (-2.01,-0.36)
          .. controls (-2.04,-0.60) and (-1.81,-0.71) .. (-1.61,-0.58)
          .. controls (-1.45,-0.47) and (-1.27,-0.53) .. (-1.23,-0.73);
        \draw[fine sulcus]
          (-2.01,-0.36) .. controls (-1.81,-0.37) and (-1.76,-0.17) .. (-1.58,-0.22);
        \draw[sulcus]
          (-2.22,-0.83) .. controls (-1.93,-0.80) and (-1.93,-1.06) .. (-1.71,-1.12)
          .. controls (-1.51,-1.16) and (-1.45,-0.97) .. (-1.24,-1.02);
        \draw[sulcus]
          (-1.91,-1.53) .. controls (-1.76,-1.35) and (-1.49,-1.43) .. (-1.48,-1.66)
          .. controls (-1.47,-1.88) and (-1.19,-1.89) .. (-1.07,-1.67)
          .. controls (-0.98,-1.51) and (-0.81,-1.49) .. (-0.69,-1.59);
        \draw[sulcus]
          (-1.40,-2.15) .. controls (-1.17,-2.09) and (-1.15,-2.32) .. (-0.95,-2.23)
          .. controls (-0.77,-2.14) and (-0.88,-1.96) .. (-0.69,-1.92);
        \draw[sulcus]
          (-0.76,-2.53) .. controls (-0.77,-2.31) and (-0.43,-2.24) .. (-0.48,-2.04)
          .. controls (-0.52,-1.89) and (-0.34,-1.72) .. (-0.23,-1.75);
        \draw[fine sulcus]
          (-0.25,2.57) .. controls (-0.46,2.51) and (-0.27,2.27) .. (-0.35,2.17);
        \draw[fine sulcus]
          (-1.73,2.27) .. controls (-1.57,2.40) and (-1.57,2.12) .. (-1.49,2.09);
        \draw[fine sulcus]
          (-2.29,1.80) .. controls (-2.06,1.89) and (-1.98,1.68) .. (-2.05,1.58);
        \draw[fine sulcus]
          (-2.38,0.24) .. controls (-2.20,0.37) and (-2.18,0.11) .. (-2.11,0.10);
        \draw[fine sulcus]
          (-2.22,-0.63) .. controls (-2.17,-0.47) and (-2.40,-0.40) .. (-2.38,-0.27);
        \draw[fine sulcus]
          (-1.93,-1.28) .. controls (-1.70,-1.21) and (-1.74,-1.40) .. (-1.63,-1.43);
        \draw[fine sulcus]
          (-1.65,-1.95) .. controls (-1.62,-1.77) and (-1.82,-1.80) .. (-1.78,-1.66);
        \path[fill=brainrim!60!brainmid]
          (-0.11,1.13)
          .. controls (-0.33,1.50) and (-0.79,1.42) .. (-1.05,1.02)
          .. controls (-1.25,0.73) and (-1.14,0.37) .. (-0.82,0.22)
          .. controls (-0.58,0.07) and (-0.56,-0.29) .. (-0.73,-0.66)
          .. controls (-0.89,-1.02) and (-1.09,-1.15) .. (-0.97,-1.47)
          .. controls (-0.62,-1.49) and (-0.34,-1.03) .. (-0.19,-0.70)
          .. controls (0.07,-0.11) and (0.01,0.67) .. (-0.11,1.13) -- cycle;
        \path[fill=white,draw=brainfold,line width=0.65pt]
          (-0.10,0.94)
          .. controls (-0.33,1.22) and (-0.64,1.20) .. (-0.84,0.93)
          .. controls (-1.03,0.69) and (-0.98,0.48) .. (-0.73,0.38)
          .. controls (-0.40,0.24) and (-0.38,-0.17) .. (-0.52,-0.52)
          .. controls (-0.63,-0.81) and (-0.82,-1.05) .. (-0.82,-1.25)
          .. controls (-0.53,-1.13) and (-0.30,-0.75) .. (-0.20,-0.48)
          .. controls (0.01,0.05) and (-0.01,0.58) .. (-0.10,0.94) -- cycle;
      \end{scope}
    \end{scope}
  }
\end{scope}
\foreach \xx in {0.65,1.20,1.75,2.30,2.85} {
  \draw[field line,incident] (\xx,2.05) -- (\xx,5.25);
}
\draw[direction,incident] (0.28,3.65) -- (3.73,3.65);
\node[above,text=incident] at (3.35,3.65) {$\mathbf d$};
\node[font=\large] at (1.76,5.80) {$u^{\mathrm i}(\mathbf{x};\mathbf d)$};
\begin{scope}[shift={(6.55,3.65)}]
  \foreach \rr in {2.95,3.49,4.03,4.57,5.11} {
    \draw[field line,scattered]
      (0:\rr) arc[start angle=0,end angle=33,radius=\rr];
  }
  \draw[direction,scattered] (16.5:2.68) -- (16.5:5.65);
  \node[above right,text=scattered] at (16.5:5.65) {$\widehat{\mathbf x}$};
\end{scope}
\node[font=\large] at (10.95,6.85) {$u^{\mathrm s}(\mathbf{x};\mathbf d)$};
\node[right,font=\large] at (7,1) {$D$};
\end{tikzpicture}}

    \caption{Schematic of the inverse scattering problem for a brain-shaped inhomogeneity (D). An incident plane wave \(u^{\mathrm i}\) propagating in direction \(\mathbf d\) impinges upon the scatterer and generates an outgoing scattered field \(u^{\mathrm s}\). }
    \label{fig:schematic}
\end{figure}

The incident  wave interacts with the scatterer and the total acoustic field $u(\bfx)$ decomposes as the sum of the incident and scattered fields. The latter field is denoted by  $u^{\mathrm{s}}$. The total field satisfies the Helmholtz equation in all $\mathbb{R}^2$ together with the aforementioned decomposition:
\begin{align}
\Delta u + \kappa^2 n^2 u&= 0 \quad \text{in } \mathbb{R}^2, \label{eqn:helmholtz} \\
    u^{\mathrm{i}} + u^{\mathrm{s}}&=u \quad \text{in } \mathbb{R}^2. \label{eqn:uius}
\end{align}
 To ensure that the scattered field is outgoing, we impose the standard Sommerfeld radiation condition
\begin{equation}
\label{eqn:sommer}
\lim_{r \to \infty} r^{1/2} 
\left( \frac{\partial u^{\mathrm{s}}}{\partial r} - i\kappa u^{\mathrm{s}} \right) = 0,
\end{equation}
where $r = \|\bfx\|$ and the limit is uniform in all unit direction vectors $\hat{\bfx} = \bfx / \|\bfx\|$. Equations \eqref{eqn:helmholtz}-\eqref{eqn:sommer} constitute the forward scattering problem and admit a unique solution~\cite{colton+kress4}.

Since the scattered field satisfies the Helmholtz equation and the Sommerfeld radiation condition, it admits the asymptotic expansion
\begin{equation}
\label{eqn:farp}
u^{\mathrm{s}}({\bfx}; \bfd) = 
\frac{\exp(i\kappa r)}{\sqrt{r}} 
\left( u_\infty(\hat{\bfx}; \bfd) + O\!\left( \frac{1}{r} \right) \right), 
\qquad r \to \infty,
\end{equation}
where $u_\infty : \mathbb{S}^1 \times \mathbb{S}^1 \to \mathbb{C}$ is the far-field pattern of the scattered wave~\cite{colton+kress4} (here $\mathbb{S}^1$ is the unit circle). The notation emphasizes that the pattern depends on both the direction of measurement $\hat{\bfx}$ and the direction of incidence $\bfd$.  The far-field pattern is the data for the 
inverse problem that can, in principle,  be measured by sensors far from the scatterer.

In Section~\ref{sec:obstacles}, we consider obstacle scattering in which the 
scatterers are impenetrable with the Neumann boundary condition on the boundary of $D$, so that the total and scattered fields are still related by (\ref{eqn:uius}) which now holds outside $D$, together with the Sommerfeld radiation condition (\ref{eqn:sommer}).  In addition
\begin{align}
    \Delta u +\kappa^2 u&=-\delta_z\quad\mbox{in }\mathbb{R}^2\setminus \overline{D},    \label{eq:helm_ext} \\ 
\frac{\partial u}{\partial\nu}&=0\quad\mbox{on }\partial D,
\label{eq:soundhard}
\end{align}
Here $\delta_z$ is the Dirac delta at $\bfz$,  $\partial D$ is the boundary of $D$ and $\nu$ is the unit outward normal to $D$.  Under suitable conditions on the boundary of $D$, this problem also has a unique solution (provided $\mathbb{R}^2\setminus \overline{D}$ is connected)~\cite{colton+kress4}.

For a given bounded penetrable scatterer $D$ and refractive index $n(\bfx)$ in the case of a penetrable scatterer, we use the finite element method implemented in NGSolve~\cite{netgen} on a triangular grid tahtb is terminated by the PML~\cite{berenger} as
implemented in NGSolve, to compute the near-field for $u$, and then extend to the far-field using the technique in~\cite{Monk+Suli:98}.  More details can be found in~\cite{DMLM26}.

For a scatterer consisting of disjoint sound-hard disks, we instead expand the scattered field in
outgoing cylindrical waves centered at each disk:
\begin{equation}
     u^s(\bfx;\bfd)\approx\sum_{j=1}^{J}\sum_{m=-M}^{M}
 b_{jm}(\bfd)H_m^{(1)}(\kappa|\bfx-\mathbf c_j|)
 e^{im\arg(\bfx-\mathbf c_j)},
\end{equation}
where $J$ is the number of disks, $\mathbf c_j$ is the center of disk
$j$, and $M$ is the truncation order.
Graf's addition theorem~\cite[Eq.~(3.111)]{colton+kress4} expresses the field
scattered by each disk in coordinates centered at the other disks.
Imposing the Neumann conditions then gives a coupled linear system
for the expansion coefficients, including multiple scattering.
Applying the large-distance asymptotics of the Hankel functions gives
\begin{equation}
     u_\infty(\hat\bfx;\bfd)\approx\sqrt{\frac{2}{\pi\kappa}}e^{-i\pi/4}
 \sum_{j=1}^{J}\sum_{m=-M}^{M}b_{jm}(\bfd)
 e^{-i\kappa\hat\bfx\cdot\mathbf c_j}e^{im(\theta-\pi/2)},
 \qquad \hat\bfx=(\cos\theta,\sin\theta).s
\end{equation}
Thus the far-field is evaluated directly from the expansion
coefficients, with the normalization in Eq.~\eqref{eqn:farp}.  For near-field scattering data from penetrable disks we
also use special function expansion inside the component discs making up $D$.

\subsection{Inverse problem}

For the model inverse time harmonic acoustic scattering problem we shall consider in most of this paper,  the multistatic far-field pattern $u_{\infty}(\hat\bfx;\bfd)$ for all $\hat\bfx,\bfd\in \mathbb{S}^1$ is assumed to be known measured data. For some upcoming experiments, we corrupt this data by random error to simulate by measurement noise as specified by \eqref{eq:measurement-noise}. From this data, it is then desired to determine a quantitative approximation to $n^2(\bfx)$. 

Quite generally, it can be shown that the refractive index is 
uniquely determined by an exact knowledge of noise-free multistatic far-field data~\cite{bukhgeim} at a single frequency.  However,  the problem is ill-posed. To see this note that the exact far-field pattern
is an analytic function of its arguments.  Thus, any non-smooth
perturbation of the exact far-field pattern, say by $L^2$ noise,
takes it outside the class of far-field patterns  and no solution exists.
To compute a solution of the inverse problem, it is therefore usual to use an optimization approach with regularization~\cite{colton+kress4}.  Such an approach is computationally time consuming because each step of the optimization algorithm (often, Gauss-Newton) requires at least the solution of a new forward problem.  An alternative approach is to linearize the forward problem about a background field resulting in the Born approximation of the 
scattered field which results in an ill-posed linear inverse problem that can then be regularized in the standard way.  This approach can be very successful provided the scatterer is not too strong in a way that we shall quantify shortly.

To develop the Born approximation, define the contrast $\mu(\bfx)=n^2(\bfx)-1$ and note that $\mu(\bfx)=0$ for $\bfx$ outside the scatterer. If $\sigma=\kappa^2 \Vert \mu\Vert_{L^{\infty}(\mathbb{R}^2)}$ is small, using the Lippmann-Schwinger equation for the field $u$, we can derive the following approximation of the scattered field, which improves as $\sigma\to 0$:
\begin{equation}
    \label{bornu}
    u^s(\textbf{x};\textbf{d})\approx \kappa^2\int_{\mathbb{R}^2}\Phi(\textbf{x},\textbf{z})\mu(\textbf{z})u^i(\textbf{z};\textbf{d})\,d\textbf{z},\quad \forall\:\textbf{x}\in \mathbb{R}^2, \forall\:\textbf{d}\in \mathbb{S}^1.
\end{equation}
Here $\Phi(\mathbf x,\mathbf z):=\frac{i}{4}H_0^{(1)}(\kappa\|\mathbf x-\mathbf z\|)$  denotes the fundamental solution of the background Helmholtz equation $\Delta_x \Phi +\kappa^2\Phi=-\delta_z$ where $\delta_z$ is the Dirac delta at $\bfz$. Using \eqref{bornu} and letting $|\bfx|\to \infty$, we obtain the \emph{Born approximation of the far-field pattern}:
\begin{equation}
    \label{born}
    u_{\infty}(\hat{\textbf{x}},\textbf{d})\approx \frac{\kappa^2}{\sqrt{8\pi \kappa} }\exp\left(\frac{\iu\pi}{4}\right)\int_{\mathbb{R}^2}\exp(\iu \kappa(\textbf{d}-\hat{\textbf{x}})\cdot \textbf{z})\mu(\textbf{z})\,d\textbf{z}\defeq (\mathscr{B}\mu)(\hat{\textbf{x}},\textbf{d}).
\end{equation} 
This defines the Born operator $\mathscr{B}:L^2(D)\to L^2(\mathbb{S}^1\times\mathbb{S}^1)$, where the two angular variables specify the observation
and incidence directions.
In this formula, the accuracy of the Born approximation improves as $\sigma$ decreases.  For small  $\sigma$ the scattering is said to be weak.

Note that since $\Vert \bfd\Vert=1$ and $\Vert \hat{\mathbf x}\Vert=1$, $\mathscr{B}$ represents a band limited two dimensional Fourier transform known only for
transform parameters with $\ell^2$ norm less than or equal
 to $2\kappa$.  Thus, ${\mathscr B}$ is compact and finding ${\mathscr B}^{-1}$ is ill-posed.

If we want to use a neural network to approximate the solution of the inverse problem, the lack of continuity of the unregularized inverse map implies that
the standard universal approximation theorem of neural networks~\cite{cybenko} does not apply. This is a point that we shall return to in the next section where we give examples of hallucinations or instability. Such artifacts can occur~\cite{siam-rev} when a neural network is applied to compute this inverse.

When using the Born approximation, it is  usual to use regularization to stabilize the inversion and  render the approximate linearized inverse map well-posed. 
To obtain a practical algorithm the data is first discretized by assuming that we have available measurements of the scattered field in discrete directions uniformly distributed around the scatterer.  Assuming there are $N_{\rm m}>0$ measurements, let \begin{equation}
\theta_j:=2\pi j/N_{\rm m}, \quad \hat\bfx_j:=(\cos\theta_j,\sin\theta_j),\quad j=1,\ldots,N_{\rm m},
\label{theta}
\end{equation}
We assume incident directions are identical to measurement directions.  We thus have available a data matrix ${U}$ with
\begin{equation}  U_{\ell,m}:=u_\infty(\hat\bfx_\ell,\hat\bfx_m),\quad 1\leq \ell,m\leq N_{\rm m}.
\end{equation}

In the upcoming noise experiments, we perturb the complex measurements
before separating their real and imaginary parts into input channels.
Writing $\bfv=\operatorname{vec}{U}\in\mathbb C^{N_{\rm m}^2}$ where $\operatorname{vec}$ vectorizes $U$ column by column,  we use
\begin{equation}
  \bfv_j^\delta:=\bfv_j(1+\delta\mathcal E_j),\qquad
  \mathcal E_j:=\frac{a_j+i b_j}{\sqrt2},\qquad
  a_j,b_j\sim\mathcal N(0,1),
  \label{eq:measurement-noise}
\end{equation}
where all $a_j$ and $b_j$ are independent. Thus
$\mathbb E|\mathcal E_j|^2=1$, and the real and imaginary
parts of $\mathcal E_j$ each have variance $1/2$. This is multiplicative noise:
$\mathbb E|v_j^\delta-v_j|^2=\delta^2|v_j|^2$.
For use with standard neural network software (in our case PyTorch) the real and imaginary parts of $\bfv^\delta$ are further vectorized to be  used as separate real channels to the relevant neural networks.

To approximate the inverse problem we assume that there is a square search region known a priori to contain the unknown scatterer.  In this paper the region is denoted $\Omega=[-1,1]^2$. So we know a priori that $n(\bfx)=1$ outside $\Omega$. The search domain $\Omega$ is covered by a uniform $N_g\times N_g$
grid of points 
\begin{equation}\label{zpq}
\bfz_{p,q}
:=\bigl(-1+(p-1)h,\,-1+(q-1)h\bigr),\quad \;1\leq p,q \leq N_{g}\end{equation}
where \(h=2/(N_g-1)\).
Applying quadrature to \eqref{born}, and noting $\mu(\bfx)=0$ outside $\Omega$, gives the discrete Born approximation
\begin{equation}
    u_{\infty}(\hat{\textbf{x}}_\ell,\hat{\textbf{x}}_m)\approx \frac{\kappa^2}{\sqrt{8\pi \kappa} }\exp\left(\frac{\iu\pi}{4}\right)h^2\sum_{p,q=1}^{N_g}\exp(\iu \kappa(\hat{\textbf{x}}_m-\hat{\textbf{x}}_\ell)\cdot \textbf{z}_{p,q})\mu(\textbf{z}_{p,q}),\quad 1\leq  \ell,m\leq N_m.
\end{equation}
After vectorization, the right hand side of the above approximation defines the Born matrix $\mathscr{B}_{a}\in \mathbb{C}^{N_m^2\times N_g^2}$.
In our studies $N_g=100$, and we use $N_{\rm m}=100$ sources and receivers.

As before, let $\bfv\in\mathbb C^{N_{\rm m}^2}$ denote the vectorized
far-field measurements and let
$\mathscr B_a\in\mathbb C^{N_{\rm m}^2\times N_g^2}$
denote the discrete Born matrix. The Born approximation is \(v\simeq\mathscr B_a\mu,\)
where $\mu\in\mathbb R^{N_g^2}$ is the vectorized contrast.
For $\gamma>0$, define the complex Tikhonov inverse
\begin{equation}
T_{a,\gamma}:=(\mathscr B_a^*\mathscr B_a+\gamma I)^{-1}\mathscr B_a^*,
\label{eq:inverse}
\end{equation}
where $*$ denotes conjugate transpose. The reconstructed
contrast image is
\[
\tilde\mu
:=\Re\bigl(T_{a,\gamma}v^\delta\bigr).
\]
The neural networks receive the two-channel real array
$\bfy^\delta=(\Re v^\delta,\Im v^\delta)
\in\mathbb R^{2N_{\rm m}^2}$.

\section{Neural Networks} \label{sec:New_NN}
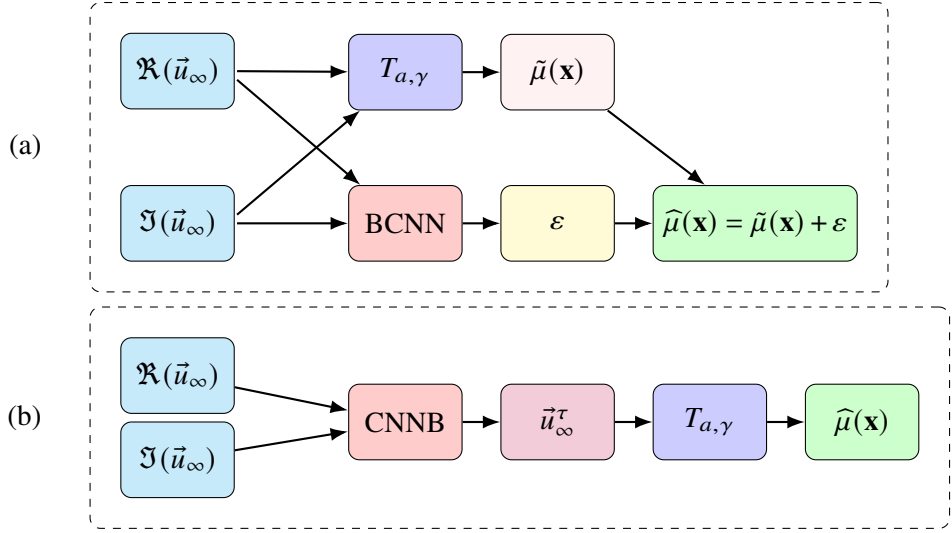
\begin{figure}[t]
    \centering
   \begin{tikzpicture}[
    node distance=1cm and 0.5cm,
    box/.style={draw, rounded corners, minimum width=1.5cm, minimum height=1cm, align=center},
    arrow/.style={-Latex, thick},
    dashedbox/.style={draw=black, dashed, rounded corners, inner sep=0.4cm}
    ]

\node[box, fill=cyan!20] (ReU2) {$\Re(\vec{u}_{\infty})$};
\node[right=-0.24cm of ReU2](RightBorderReU2){};
\node[box, fill=cyan!20, below=1cm of ReU2] (ImU2) {$\Im(\vec{u}_{\infty})$};
\node[right=-0.24cm of ImU2](RightBorderImU2){};
\node[box, fill=blue!20, right=1.5cm of ReU2.east|-ImU2.east, yshift=2cm] (B2) {$T_{a,\gamma}$};
\node[box, fill=pink!20, right=of B2] (TildeM) {$\tilde{\mu}(\textbf{x})$};
\node[box, fill=red!20, below right=of ReU2, xshift=1cm] (NN2) {BCNN};
\node[box, fill=yellow!20, right=of NN2] (Epsilon) {$\varepsilon$};
\node[box, fill=green!20, right=of Epsilon] (M2) {$\widehat\mu(\textbf{x})=\tilde{\mu}(\textbf{x})+\varepsilon$};

\draw[arrow] (RightBorderReU2) -- (B2);
\draw[arrow] (RightBorderImU2) -- (B2);
\draw[arrow] (RightBorderReU2) -- (NN2);
\draw[arrow] (RightBorderImU2) -- (NN2);
\draw[arrow] (B2) -- (TildeM);
\draw[arrow] (TildeM) -- (M2);
\draw[arrow] (NN2) -- (Epsilon);
\draw[arrow] (Epsilon) -- (M2);

\node[dashedbox, fit=(ReU2) (ImU2) (B2) (TildeM) (NN2) (Epsilon) (M2)] (box2) {};
\node[left=0.5cm of box2.west] {(a)};

\node[box, fill=cyan!20, below=1cm of ImU2] (ReU1) {$\Re(\vec{u}_{\infty})$};
\node[box, fill=cyan!20, below=0.1cm of ReU1] (ImU1) {$\Im(\vec{u}_{\infty})$};
\node[box, fill=red!20, right=1.5cm of ReU1.east|-ImU1.east, yshift=0.5cm] (NN1) {CNNB};
\node[box, fill=purple!20, right=of NN1] (VF1) {$\vec{u}_{\infty}^{\tau}$};
\node[box, fill=blue!20, right=of VF1] (B1) {$T_{a,\gamma}$};
\node[box, fill=green!20, right=of B1] (M1) {$\widehat\mu(\textbf{x})$};

\draw[arrow] (ReU1) -- (NN1);
\draw[arrow] (ImU1) -- (NN1);
\draw[arrow] (NN1) -- (VF1);
\draw[arrow] (VF1) -- (B1);
\draw[arrow] (B1) -- (M1);

\node[dashedbox, fit=(ReU1) (ImU1) (NN1) (VF1) (B1) (M1)] (box3) {};
\node[left=0.5cm of box3.west] {(b)};

\end{tikzpicture}
\caption{Illustrations of the networks BCNN and CNNB adapted from \cite{DMLM26} (with permission), using the present contrast notation. The real part is taken after applying $T_{a,\gamma}$. The top panel illustrates BCNN in which a correction is added
after applying the inverse Born approximation.  The bottom panel illustrates 
CNNB in  which the far-field data is modified before applying the inverse Born approximation. }
\label{oldnet}
\end{figure}

Generally, as $\sigma$ increases and we move out of the weak scattering regime,
the accuracy of the regularized inverse Born approximation deteriorates.  Therefore,
in \cite{DMLM26}, we proposed using neural networks to correct the inverse Born approximation.  For the correction networks we chose convolutional neural networks (CNNs). In particular, we proposed either to try to train the neural network to identify the Born approximation error in a post-processing step by adding a learned image correction,
or pre-process the far-field data using a neural network before applying the regularized Born approximation.  The former case is called BCNN and a graphic
of the network is shown in Fig.~\ref{oldnet})a).  The second, pre-processing version is called CNNB and is shown in the bottom panel of Fig.~\ref{oldnet}.    Both correction networks receive the
far-field scattering data as input~\cite{DMLM26}.  Details of determining the network parameters such as number of layers and number of neurons per layer can also be found in that paper, as well as a discussion of training.

For all the numerical experiments reported here, we use the wave number $\kappa=16$ giving a wavelength of 0.39.  The BCNN and CNNB neural networks were trained on 16,000 far-field patterns (without additional noise beyond finite element error) with a further 4,000 far-field patterns used for validation.  The training, validation and test scatterers were taken to be one, two or three nonintersecting disks (the number chosen uniformly at random for each training sample). The radius of each disc was drawn from ${\cal U}(0.1,0.3)$ and position from  ${\cal U}(-.7,.7)^2$.  

The squared refractive index $\eta=n^2$ is taken to be a real constant
on each disc and chosen from ${\cal U}(1.1,1.8)$; the corresponding
contrast is $\mu=\eta-1$. For the wave number and scatterers in this study, the upper end of the range for $\eta$ is not weak scattering as we shall see in the numerical results.

The BCNN and CNNB models use GELU activations and
squared-$\ell^2$ training losses. Their architectures were selected
using validation data, with a separate set of 4,000 samples used
for testing; see~\cite{DMLM26} for details.

In \cite{DMLM26} we compared the two corrected Born networks to a pure CNN network as well as to a state of the art Neural Inverse Operator (NIO) network~\cite{NIO} and a Born inspired special purpose network BCR-Net~\cite{Fan:22} trained on our data.  Generally, we found that our CNN based correction algorithms performed similarly to NIO and BCR-Net (although with less sharp boundaries for the scatterers) for in-distribution scattering. For noisy data, CNNB was the most stable for higher noise levels and both BCR-Net and NIO showed some hallucinations, while the pure CNN network was less stable to higher noise levels than the correction networks.  The major difference was seen when we applied the neural networks to some more complex out of distribution scatterers. In these cases the CNN based methods produced higher-fidelity inversions than NIO or BCR-Net.
We shall compare the new networks in this paper also to NIO and BCR-Net.

A primary goal of this paper is to improve the accuracy of the
Born based networks and obtain sharper images that are closer to, or even improve on, those from NIO or BCR-Net.
Our networks CNNB and BCNN target modifying the scattering data and estimating error respectively.  It is natural to
extend this paradigm to both pre-process the scattering data and correct the image with more advanced networks. Our first network, 
called \emph{BornNet},  combines a U-Net that corrects the
far-field data, a fixed regularized inverse Born approximation, and a residual network that refines
the resulting image. The U-Net includes self-attention layers.
We denote the data and image networks by $C_\theta$ and $Q_\psi$,
respectively, where $\theta$ and $\psi$ are their learned parameters.
With $x=(\Re v,\Im v)$ denoting the two real input channels, the
reconstruction is encoded by $\hat{\mu}$ defined by
\begin{equation}
  z=\Re\bigl(T_{a,\gamma}(C_\theta(x))_c\bigr),\qquad
  \widehat\mu=Q_\psi(z).
  \label{eq:bornnet}
\end{equation}

\begin{figure}[t]
\centering
\begin{tikzpicture}[
    >=Stealth,
    tensor/.style={draw=black!60, fill=blue!5, rounded corners=3pt,
        align=center, text width=2.7cm, minimum height=0.7cm,
        inner sep=5pt, font=\footnotesize},
    module/.style={tensor, fill=green!10, minimum height=0.85cm},
    physics/.style={tensor, fill=red!10, text width=2.9cm,
        minimum height=1.25cm},
    learned/.style={physics, fill=orange!10},
    arrow/.style={->, thick},
    lossarrow/.style={->, dashed, draw=black!60}
]
\node[font=\small\bfseries, align=center] at (-3.9,1.0)
    {(a) BornNet and LinearBornNet};
\node[tensor] (inputA) at (-3.9,0) {Scattering data $x$\\two real channels};
\node[module] (unet) at (-3.9,-1.2) {Data correction\\U-Net $C_\theta$};
\node[physics] (fixed) at (-5.7,-3.0)
    {BornNet\\fixed Born inverse\\$z=\Re(T_{a,\gamma}c_c)$};
\node[learned] (linear) at (-2.1,-3.0)
    {LinearBornNet\\learned linear map\\$z=W_2W_1c_f$};
\node[tensor] (intermediate) at (-3.9,-4.65)
    {Intermediate contrast $z$\\$100\times100$};
\node[module] (refiner) at (-3.9,-5.9)
    {Image refinement $Q_\psi$\\3 residual blocks};
\node[tensor] (imageA) at (-3.9,-7.15)
    {Predicted contrast $\widehat\mu$};

\draw[arrow] (inputA) -- (unet);
\coordinate (split) at (-3.9,-2.05);
\draw (unet.south) -- (split)
    node[midway,right,font=\footnotesize] {$c=C_\theta(x)$};
\draw[arrow] (split) -| (fixed.north);
\draw[arrow] (split) -| (linear.north);
\coordinate (merge) at (-3.9,-4.0);
\draw (fixed.south) |- (merge);
\draw (linear.south) |- (merge);
\draw[arrow] (merge) -- (intermediate.north);
\draw[arrow] (intermediate) -- (refiner);
\draw[arrow] (refiner) -- (imageA);

\node[font=\small\bfseries] at (3.0,1.0) {(b) BornAE};
\node[tensor] (inputB) at (3.0,0) {Scattering data $x$\\two real channels};
\node[module, text width=3.2cm] (encoder) at (3.0,-1.2)
    {LinearBornNet encoder\\from panel (a)};
\node[tensor] (imageB) at (3.0,-2.4)
    {Predicted contrast $\widehat\mu$};
\node[module, fill=red!10, text width=3.2cm] (forward) at (3.0,-3.65)
    {Fixed forward Born map\\$\mathscr B_a\widehat\mu$};
\node[module, text width=3.2cm] (corrector) at (3.0,-4.9)
    {Field correction\\2 residual blocks\\+ output layer};
\node[tensor] (fieldB) at (3.0,-6.15)
    {Reconstructed data\\$\widehat x=D_\phi(\widehat\mu)$};
\node[tensor, fill=gray!10, text width=4.2cm, below=0.7cm of fieldB] (loss)
    {Training loss $\mathcal{L}_2$\\[3pt]};

\draw[arrow] (inputB) -- (encoder);
\draw[arrow] (encoder) -- (imageB);
\draw[arrow] (imageB) -- (forward);
\draw[arrow] (forward) -- (corrector);
\draw[arrow] (corrector) -- (fieldB);
\draw[lossarrow] (fieldB) -- (loss);
\draw[lossarrow] (imageB.east) -- ++(2.1,0) |- (loss.east);
\begin{scope}[on background layer]
\node[draw=black!50, dashed, rounded corners, inner sep=8pt,
    fit=(forward)(corrector)(fieldB),
    label={[font=\scriptsize]left:training-only decoder}] {};
\end{scope}
\end{tikzpicture}
\caption{The three new networks. (a) BornNet and LinearBornNet differ
in the map from corrected data to the intermediate contrast.
Here $c_c$ combines the two real channels of $c$ into a complex
vector, while $c_f$ stacks them into one real vector. The learned
map has dimensions $20{,}000\to1{,}024\to10{,}000$.
Both $z$ and $\widehat\mu$ are supervised during training.
(b) BornAE uses the LinearBornNet encoder, with the forward Born map
and learned field correction used only during training. Only
$\widehat\mu$ is needed at inference. The diagram shows the far-field
operators; the near-field scaling is given in Section~\ref{sec:near}.}
\label{fig:new}
\end{figure}
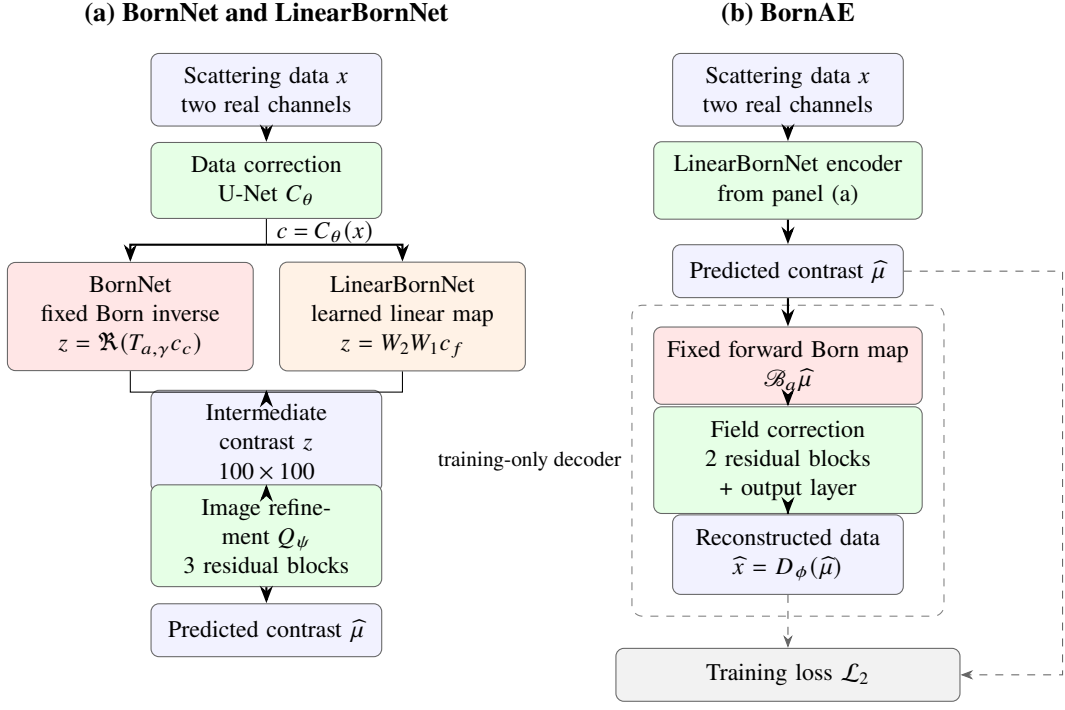
Here the subscript $c$ combines the two real output channels into
complex data, and the real part of the fixed inverse Born output
gives a contrast image. The intermediate image $z$ and the final reconstruction
are both compared with the true contrast in the squared-error training
loss so that
\begin{equation}
    {\cal L}_1= \frac{1}{N_{g}^2}\lVert\widehat\mu-\mu\rVert_2^2
 +\frac{1}{N_{g}^2}\lVert z-\mu\rVert_2^2.
\end{equation}
The two networks in \emph{BornNet} are trained jointly.

The next network is intended to investigate the question: Is the Born approximation needed? The motivation here is that for near-field imaging, particularly in 3D, the map $\mathscr{B}$ is no longer a band-limited Fourier transform, and so it is more difficult to regularize and invert the Born approximation. The new network, called \emph{LinearBornNet}, uses the same two pre- and post-processing networks as BornNet,
but replaces the fixed inverse layer by the product $W_2W_1$ of two learned real matrices. The first, $W_1$, maps
the 20,000 corrected-data values output by the U-Net to 1,024 coefficients; the
second, $W_2$, maps those coefficients to the 10,000 image pixels.  There is no
activation between the two matrices, so their product has rank at most
1,024; see Fig.~\ref{fig:new}(a).

No effort was made to determine the ideal 
number of rows for $W_1$.  The loss still includes the intermediate-image error.
This use of a low rank approximation is also motivated by regularization by truncation that can also be used to stabilize the Born approximation.

The final network considered in this paper is called \emph{BornAE} and is designed to explicitly incorporate the Born approximation. It uses
the forward Born approximation as a consistency regularizer.  It uses the LinearBornNet architecture as its image encoder and
adds a field decoder.  The decoder applies the fixed forward Born matrix to the
predicted contrast, then corrects the resulting field with two residual
blocks and a two-channel convolutional output.  See Fig~\ref{fig:new}(b).  
This network avoids computing $T_{a,\gamma}$, but the regularization afforded by the decoder is  weak.   If $D_\phi$ denotes
this decoder and $N_g^2$ is the number of image pixels, its training loss is
\begin{equation}
 \mathcal L_2=
 \frac{1}{N_g^2}\lVert\widehat\mu-\mu\rVert_2^2
 +\frac{1}{2m}\lVert D_\phi(\widehat\mu)-x\rVert_2^2
 +\frac{10^{-4}}{N_g^2}\lVert\widehat\mu\rVert_1.
 \label{eq:bornae-loss}
\end{equation}
The first two terms measure image error and field-reconstruction error;
the last is an $\ell^1$ penalty on the reconstructed contrast.
Here $m=N_{\rm m}^2$, and the far-field measurements and decoder
outputs are un-normalized.
Unlike LinearBornNet, BornAE does not penalize the intermediate-image
error. The field loss uses the decoder’s corrected Born prediction. At inference, only the encoder is needed to reconstruct the image.

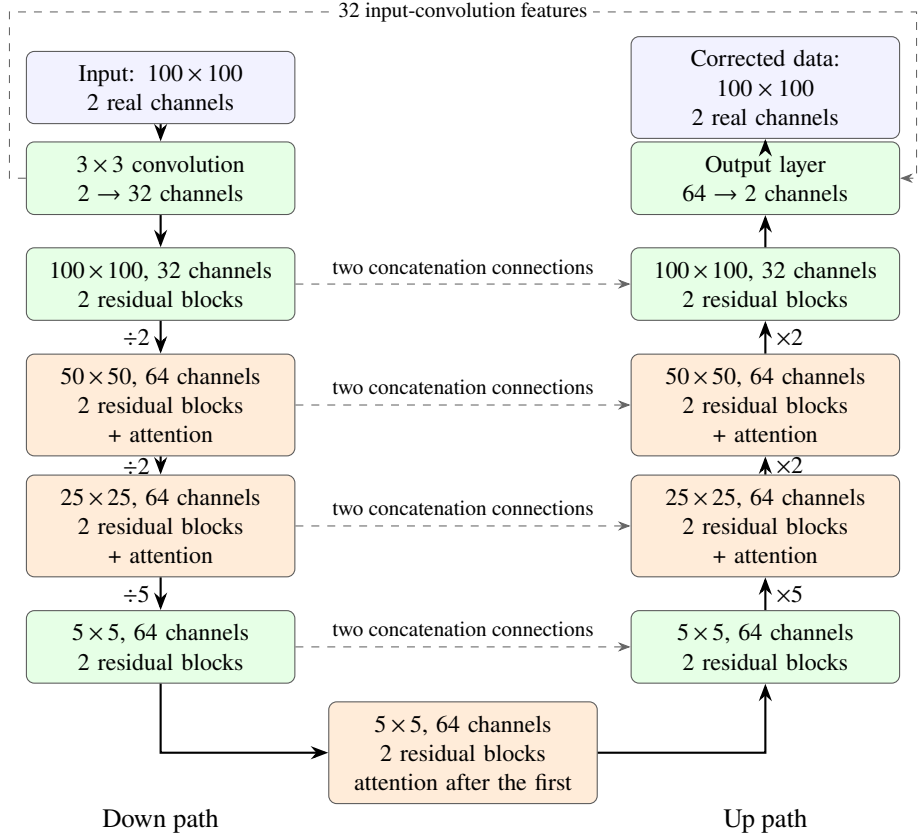
\begin{figure}[t]
\centering
\begin{tikzpicture}[
    >=Stealth,
    block/.style={draw=black!65, rounded corners=3pt, fill=green!10,
        text width=3.2cm, minimum height=0.9cm, inner sep=5pt,
        align=center, font=\footnotesize},
    attn/.style={block, fill=orange!15},
    io/.style={block, fill=blue!5, minimum height=0.7cm},
    arrow/.style={->, thick},
    skip/.style={->, dashed, draw=black!60}
]
\node[io] (input) at (-4,2.0) {Input: $100\times100$\\2 real channels};
\node[block] (stem) at (-4,0.8) {$3\times3$ convolution\\$2\to32$ channels};
\node[block] (d100) at (-4,-0.6) {$100\times100$, 32 channels\\2 residual blocks};
\node[attn] (d50) at (-4,-2.2) {$50\times50$, 64 channels\\2 residual blocks + attention};
\node[attn] (d25) at (-4,-3.8) {$25\times25$, 64 channels\\2 residual blocks + attention};
\node[block] (d5) at (-4,-5.4) {$5\times5$, 64 channels\\2 residual blocks};
\node[attn] (centre) at (0,-6.8) {$5\times5$, 64 channels\\2 residual blocks\\attention after the first};
\node[block] (u5) at (4,-5.4) {$5\times5$, 64 channels\\2 residual blocks};
\node[attn] (u25) at (4,-3.8) {$25\times25$, 64 channels\\2 residual blocks + attention};
\node[attn] (u50) at (4,-2.2) {$50\times50$, 64 channels\\2 residual blocks + attention};
\node[block] (u100) at (4,-0.6) {$100\times100$, 32 channels\\2 residual blocks};
\node[block] (head) at (4,0.8) {Output layer\\$64\to2$ channels};
\node[io] (output) at (4,2.0) {Corrected data: $100\times100$\\2 real channels};

\draw[arrow] (input) -- (stem);
\draw[arrow] (stem) -- (d100);
\draw[arrow] (d100) -- (d50) node[midway,left,font=\footnotesize] {$\div2$};
\draw[arrow] (d50) -- (d25) node[midway,left,font=\footnotesize] {$\div2$};
\draw[arrow] (d25) -- (d5) node[midway,left,font=\footnotesize] {$\div5$};
\draw[arrow] (d5.south) |- (centre.west);
\draw[arrow] (centre.east) -| (u5.south);
\draw[arrow] (u5) -- (u25) node[midway,right,font=\footnotesize] {$\times5$};
\draw[arrow] (u25) -- (u50) node[midway,right,font=\footnotesize] {$\times2$};
\draw[arrow] (u50) -- (u100) node[midway,right,font=\footnotesize] {$\times2$};
\draw[arrow] (u100) -- (head);
\draw[arrow] (head) -- (output);
\foreach \level in {100,50,25,5} {
    \draw[skip] (d\level.east) -- (u\level.west)
        node[midway,above,font=\scriptsize] {two concatenation connections};
}
\draw[skip] (stem.west) -- (-6.0,0.8) -- (-6.0,3.0)
    -- (6.0,3.0) -- (6.0,0.8) -- (head.east);
\node[font=\scriptsize, fill=white, inner sep=2pt] at (0,3.0)
    {32 input-convolution features};
\node[font=\small] at (-4,-7.7) {Down path};
\node[font=\small] at (4,-7.7) {Up path};
\end{tikzpicture}
\caption{The U-Net $C_\theta$ shared by all three new networks.
The two measurement arrays are processed at successively coarser
resolutions before being brought back to their original size.
On the up path, the blocks also receive the saved down-path outputs
through the dashed connections. At each resolution, these connections
pair the two blocks in reverse order and concatenate their feature
arrays (the boxes give channel counts before concatenation).
The output layer concatenates the
32 input-convolution features with the 32 up-path features, then
applies normalization, SiLU, and a $3\times3$ convolution.
The orange boxes indicate where self-attention is applied, following
both blocks at the $50\times50$ and $25\times25$ levels and the first
block at the center.
The fixed or learned inverse and image refiner follow this U-Net
and are shown separately in Fig.~\ref{fig:new}.}
\label{fig:unet}
\end{figure}

All three networks use the U-Net in Fig.~\ref{fig:unet} and the image
refiner in Table~\ref{tab:new-network-architecture}. Layer and parameter
counts are given in Table~\ref{tab:sizes}.
We represent the real and imaginary measurements by two $100\times100$
arrays, indexed by observation and incident directions. The first
convolutional layer applies 32 learned local filters to these two
arrays, producing 32 new arrays of the same size, which we call
feature \emph{channels}. Subsequent layers combine these channels
using further learned filters.

The four U-Net levels have array sizes $100\times100$, $50\times50$,
$25\times25$, and $5\times5$, with 32, 64, 64, and 64 channels,
respectively. The convolutional layers reduce the array size along
the down path, and transposed convolutions then increase it along
the up path until it reaches the original resolution. At each level, skip connections concatenate the down-path block outputs with the inputs to the corresponding up-path blocks.

We use two \emph{residual blocks} at each level on both paths and
another two at the center of the U-Net. Each residual block acts on
its input $f$ as
\[
 f\longmapsto \frac{P(f)+F(f)}{\sqrt2},
\]
where $P$ is the identity unless a learned $1\times1$ convolution
with bias is needed to match channel counts. The map $F$ contains two $3\times3$
convolutions, each preceded by group normalization and the pointwise
nonlinearity $\operatorname{SiLU}(t)=t/(1+e^{-t})$. Group normalization
centers and scales groups of feature values within each sample.
We use 32 groups, except in the first image-refinement and field-decoder
blocks, which use one.
During training, dropout randomly suppresses 10\% of the activations
before the second convolution.

In a self-attention block, each grid location receives a weighted sum
of features from all locations. The weights depend on the current
features and are computed separately by four attention heads, each
using 16 channels. The block uses $1\times1$ projections for queries,
keys, and values, and a fourth $1\times1$ projection is applied to the
combined head outputs. This projected output is added to the block
input, and the sum is divided by $\sqrt2$.
Attention follows each residual block at the
$50\times50$ and $25\times25$ levels on both paths, and the first
central block, giving nine attention blocks in total.

For a convolution, the \emph{kernel size} $K$ is the width of its
$K\times K$ learned stencil, the \emph{stride} $S$ is the number of
grid entries by which that stencil moves between output locations,
and the \emph{padding} $P$ is the number of zero entries added at
each edge before applying the stencil. An input of width $N$ then
has output width
\[
 N_{\rm out}:=\left\lfloor\frac{N+2P-K}{S}\right\rfloor+1.
\]
The first two downsampling layers use $(K,S,P)=(3,2,1)$, and the third
uses $(6,5,2)$, reducing the array width from 100 to 50, then to 25 and
finally to 5. On the up path, transposed convolutions distribute
coarse-grid features onto finer grids. Their kernel sizes and strides
are both $5,2,2$, and with zero padding they increase the array width
from 5 to 25, then to 50 and finally to 100. Their weights are learned
independently of the downsampling weights.

We pass the two corrected measurement channels from the U-Net through
the inverse layer and refine the resulting contrast image using three
residual blocks. BornAE uses two further residual blocks to
correct the predicted field during training. There is no activation
between the two learned linear maps or after the final outputs.
For near-field data, we use the same architecture with the Born
operators and field scaling described in Section~\ref{sec:near}.

\begin{table}[!htbp]
\caption{Architecture of the shared data and image networks and the
model-specific inverse and decoder stages. Channel counts give the
number of feature arrays, and unit counts give vector dimensions.
U-Net widths are block outputs, before concatenating features from
the earlier path. BornAE uses only its encoder for reconstruction.}
\centering
\small
\setlength{\tabcolsep}{4pt}
\begin{tabular}{p{0.22\linewidth}p{0.34\linewidth}p{0.35\linewidth}}
\toprule
Stage & Layers & Array size / channels or units\\
\midrule
Input & $3\times3$ convolution & $100^2$ / $2\to32$\\
U-Net down path & 8 residual blocks; 3 downsamplers &
$100^2/32$, $50^2/64$, $25^2/64$, $5^2/64$\\
U-Net center & 2 residual blocks & $5^2/64$\\
U-Net up path & 8 residual blocks; 3 upsamplers &
$5^2/64$, $25^2/64$, $50^2/64$, $100^2/32$\\
Field head & Normalization, SiLU, $3\times3$ convolution & $100^2$ / $64\to2$\\
BornNet inverse & Fixed regularized Born inverse & $20{,}000\to10{,}000$\\
Learned inverse & Two bias-free linear layers (LinearBornNet and BornAE) &
$20{,}000\to1{,}024\to10{,}000$\\
Image refiner & 3 residual blocks; normalization, SiLU, $3\times3$ head &
$100^2$ / $1\to32\to32\to32\to1$\\
BornAE decoder & Fixed Born forward map; 2 residual blocks;
normalization, SiLU, $3\times3$ head &
$100^2$ / $1\xrightarrow{\mathscr B_a}2\to32\to32\to2$\\
\bottomrule
\end{tabular}
\medskip
\label{tab:new-network-architecture}
\end{table}

\begin{table}
\caption{Network sizes. Layer counts cover convolutions, transposed
convolutions, and dense layers on all branches, including shortcut
and attention projections. Normalization layers and fixed Born maps
are excluded from the layer counts. The parameter totals include
learned normalization parameters, and the BornAE totals also include its
training-only field decoder.}
\label{tab:sizes}
\centering
\begin{tabular}{lrrr}
\toprule-
Model & Residual blocks & Learned layers & Trainable parameters\\
\midrule
BornNet & 21 & 96 & 1,936,453\\
LinearBornNet & 21 & 98 & 32,656,453\\
BornAE (training) & 23 & 104 & 32,685,739\\
\bottomrule
\end{tabular}
\end{table}

\section{Performance comparisons with far-field data}\label{sec:perf}
\subsection{Reconstruction from noisy data}\label{sec:noise}
In this section we consider reconstruction of in-distribution scatterers, but using data that has been corrupted by an increasing intensity $(\delta)$ of random noise.  Note that all training used far-field data computed by the finite element method without additional random noise.  In Table~\ref{tab:noise}, we measure the relative error in the contrast,
\begin{equation}
  E_{\mu,p}:=\frac{\lVert\widehat\mu-\mu\rVert_p}{\lVert\mu\rVert_p},\quad p=1,2
  \label{eq:contrast-error}
\end{equation}
and report their mean over the test objects.
\begin{table}[!htbp]
\caption{Noise experiments: mean relative contrast error $100E_{\mu,2}$
on 4,000 test objects (lower is better).
The first six rows are reproduced from~\cite{DMLM26}; limitations of
this comparison are discussed in the text. BornNet, LinearBornNet and
BornAE are the new networks described here. All networks were trained
without added noise. Bold entries mark the smallest mean error at
each noise level.}
\centering
\small
\begin{tabular}{lrrrrr}
\toprule
Method & $\delta=0$ & $\delta=0.1$ & $\delta=0.25$ & $\delta=0.5$ & $\delta=1$ \\
\midrule
Born ($\gamma=1$) & 67.04 & 67.04 & 67.04 & 67.04 & 67.04 \\
CNN & 35.01 & 36.10 & 41.96 & 62.67 & 134.24 \\
NIO & 32.70 & 36.88 & 63.11 & 103.49 & 148.96 \\
BCR-Net & 33.68 & 33.86 & 34.81 & 37.87 & 47.35 \\
BCNN & 36.18 & 38.41 & 48.45 & 89.34 & 380.37 \\
CNNB & 35.31 & 35.88 & 39.29 & 49.25 & 73.58 \\
BornNet & 17.75 & 19.85 & 28.38 & 55.40 & 98.49 \\
LinearBornNet & \textbf{17.57} & \textbf{17.43} & \textbf{24.11} & 39.61 & 67.01 \\
BornAE & 27.35 & 27.39 & 27.65 & \textbf{28.67} & \textbf{32.58} \\
\bottomrule
\end{tabular}
\label{tab:noise}
\end{table}
 
Starting with noise-free data ($\delta=0$ in (\ref{eq:measurement-noise})), the first data column of Table 
\ref{tab:noise} suggests that the new networks indeed result in more accurate reconstructions on average. In particular BornNet and LinearBornNet have much lower error than the other networks, while BornAE decreases the error but less substantially.

\begin{figure}[!htbp]
\centering
\includegraphics[width=\linewidth]{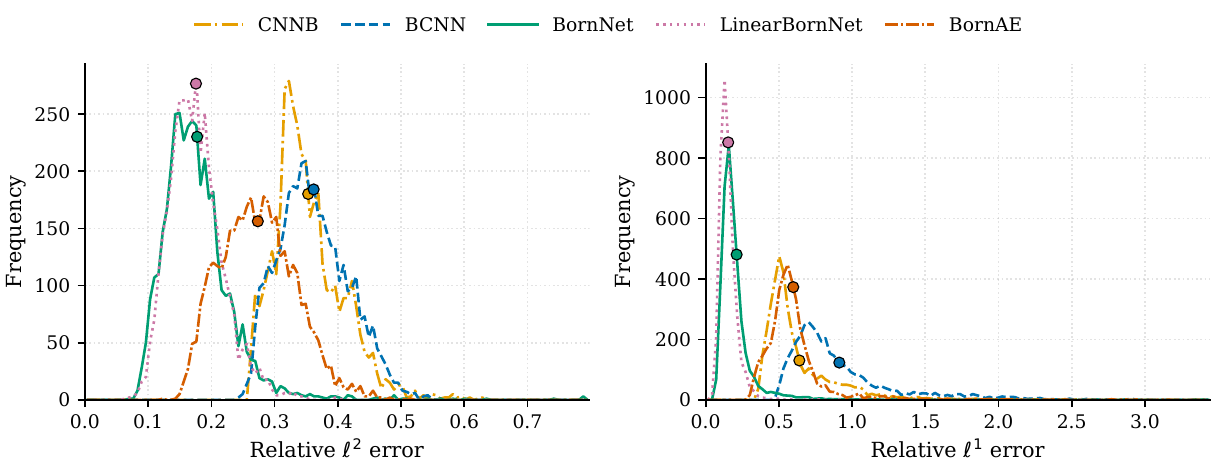}
\caption{Relative contrast errors in the $\ell^2$ (left) and $\ell^1$
(right) norms for 4,000 noiseless in-distribution test objects.  Curves show histogram
counts using the same bins for all methods within each panel;
circles mark the means.}
\label{fig:farfield-errors}
\end{figure}

The error distributions for the new networks, CNNB, and BCNN are compared in Fig.~\ref{fig:farfield-errors}.  Ideally, in this figure, the networks should show a sharp peak at small errors and a small tail of larger errors to the right.  This is indeed observed with BornNet and LinearBornNet giving very similar error behavior with average error smaller than the other networks.  BornAE also outperforms our original networks but has a wider error distribution than the other two new networks.

\begin{figure}[!htbp]
\centering
\includegraphics[width=\linewidth]{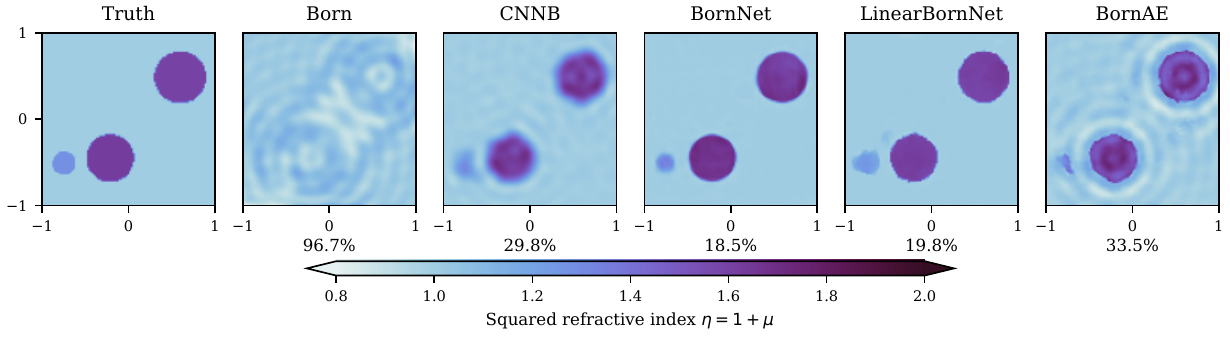}
\caption{Reconstructions of three penetrable disks from noiseless
data.  Each panel shows the squared refractive index on $[-1,1]^2$;
numbers below the reconstructions are percentage contrast errors.}
\label{fig:farfield-clean}
\end{figure}

Figure~\ref{fig:farfield-clean} shows the three-circle example from
\cite[Fig.~5]{DMLM26} using  in-distribution objects and  noise-free data.  Here and in later figures of reconstructions
we show $\eta=1+\mu$ on the common color scale
$[0.8,2]$ used in~\cite{DMLM26}; values outside this range take
the endpoint colors. Numbers
below the reconstructions give $100E_{\mu,2}$.  The Born approximation result shows very poor approximation of the two larger discs which are not weak scatterers (the small scatterer is lost completely).  It is exactly the poor performance we wish to correct.

From Figure~\ref{fig:farfield-clean} it is also clear that the BornNet
and LinearBornNet provide sharper estimates of $\mu$ compared to our previous network CNNB. All three networks indicate the presence of the small low contrast scatterer.   BornAE shows some instability in the background, but reconstructs the strong scatterers well.

For the three new networks in Table~\ref{tab:noise}, the noise-free
finite element data are corrupted using $\delta>0$ in
Eq.~\eqref{eq:measurement-noise}. The first six rows are inherited
from~\cite{DMLM26} and were not recomputed with the evaluation code
used here.
The networks were not retrained for the noisy inputs, and remain the same for all choices of $\delta$.
As we have seen before, CNNB is less affected by noise than BCNN
(Table~\ref{tab:noise}).  At $\delta=0.25$, their errors are 39.29\% and
48.45\%; at $\delta=1$, their errors are 73.58\% and 380.37\%.
This suggests that the position of the learned correction matters,
although the networks also differ in their fitted parameters and
training objectives.

BornNet and LinearBornNet improve accuracy at the lower noise levels.
At $\delta=0.25$, their mean errors are 28.38\% and 24.11\%, respectively,
both below BCR-Net's 34.81\%.  LinearBornNet has a smaller error than the
fixed-inverse BornNet at every noise level tested.
At the largest noise level, LinearBornNet is close to the regularized
Born reconstruction, while BornNet has a larger error than both it and
CNNB. BornAE trades clean accuracy for a much smaller increase in error with
noise.  Its mean error rises from 27.35\% without noise to 32.58\% at
$\delta=1$.  It has the smallest error among the reported methods at
$\delta=0.5$ and $\delta=1$, below BCR-Net's 37.87\% and 47.35\%.
The clean-data disadvantage also appears on a separate test set of
100 rectangular objects: BornAE's error is 40.07\%, compared with
33.83\% for BornNet and 34.92\% for LinearBornNet.

\begin{figure}[!htbp]
\centering
\begin{minipage}{0.32\linewidth}
\centering
\includegraphics[width=\linewidth]{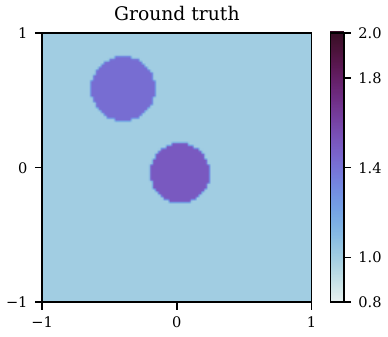}
\par\small (a) True squared refractive index
\end{minipage}
\par\medskip
\includegraphics[width=\linewidth,height=0.62\textheight,keepaspectratio]
{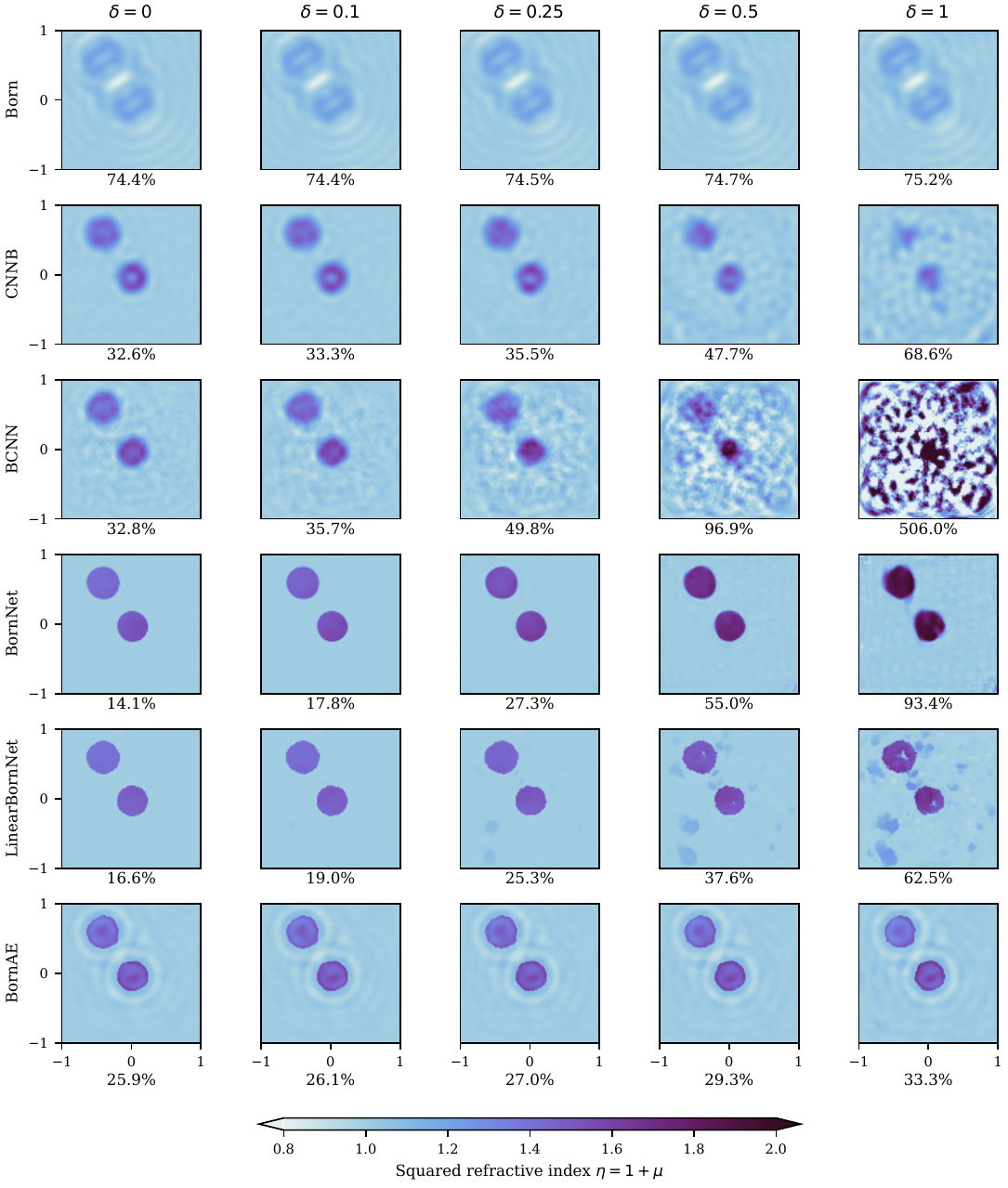}
\par\small (b) Reconstructions at different noise levels
\caption{Reconstructions of a two-circle scatterer with different
levels of noise added via Eq.~\eqref{eq:measurement-noise}.
(a) The true squared refractive indices are approximately 1.4163
and 1.4950 inside the disks, with background 1.
(b) All methods receive the same
perturbed data.  The clean (noise free) reconstructions are shown in the first
column for comparison.  Colors show squared refractive index; numbers give $100E_{\mu,2}$.}
\label{fig:farfield-noise}
\end{figure}

These comparisons use one trained network of each type, with different
training objectives and numbers of training steps.  The BornNet,
LinearBornNet and BornAE evaluations also use independent noise
realizations.  The results therefore do not isolate the effects of
attention, the inverse layer, the image network or BornAE's decoder. Detailed comparisons of the reconstructions are shown in Figure~\ref{fig:farfield-noise} which uses the two-circle object
from~\cite[Fig.~7]{DMLM26}.  For this figure, we generate a new noise
array and use it for all methods, scaling it by $\delta$.
At the largest noise level ($\delta=1$), BCNN produces large errors throughout the image suggesting instability.
BornNet still locates both disks, but their reconstructed contrasts
are too large.  LinearBornNet shows some hallucinations. BornAE has a larger clean ($\delta=0$) error but a smaller error
at the highest noise level: 25.9\% and 33.3\%, respectively,
compared with 14.1\% and 93.4\% for BornNet.

\subsection{Performance on out-of-distribution data}
One of the surprising results from \cite{DMLM26} was that the CNNB and BCNN networks, trained only on scattering from a few discs, could successfully reconstruct several out-of-distribution scatterers with apparently different structure.  
In this section, we test the new networks on these scatterers, as well as another higher contrast U-shaped scatterer, using the same protocol as in \cite{DMLM26}. In particular, for each scatterer we generate 100 realizations with
random type-specific variations. The performance of the various reconstruction algorithms reported
in Table~\ref{tab:errors} is computed from this data.

\begin{table}[h]
\caption{Average relative $\ell^2$ ($E_{\mu,2}$) and $\ell^1$ ($E_{\mu,1}$) errors (\%) in the
reconstruction of the contrast $\mu$ on 100 randomized samples of the complex OOD scatterers. Entries are reported as
mean $\pm$ standard deviation. The best average error for each problem
is in boldface.}
\centering
\renewcommand{\arraystretch}{1.2}
\begin{tabular}{|l|c c||c c|}
\hline
Class & \multicolumn{2}{c||}{\textbf{U}} & \multicolumn{2}{c|}{\textbf{Ring}} \\
\hline
Method & $\ell^2$ (\%) & $\ell^1$ (\%) & $\ell^2$ (\%) & $\ell^1$ (\%) \\
\hline
Born ($\gamma=1$) & 69.66 $\pm$ 19.80 & 119.41 $\pm$ 29.91 & 109.32 $\pm$ 11.37 & 116.41 $\pm$ 10.64 \\
NIO               & 77.87 $\pm$  8.64 & 109.20 $\pm$ 23.76 & 94.39 $\pm$ 32.38  & 98.33 $\pm$ 35.39 \\
BCR-Net           & 80.32 $\pm$  9.18 & 105.61 $\pm$ 19.88 & 73.96 $\pm$ 11.58  & 74.05 $\pm$ 16.86 \\
CNN               & 55.86 $\pm$ 10.09 & 110.55 $\pm$ 32.51 & 64.77 $\pm$ 23.97  & 65.27 $\pm$ 24.36 \\
CNNB              & 53.94 $\pm$ 9.95 & 98.37 $\pm$ 32.14 & \textbf{60.57 $\pm$ 25.32} & \textbf{63.82 $\pm$ 26.09} \\
BCNN              & 53.54 $\pm$  9.42 & 130.51 $\pm$ 35.91 & 74.55 $\pm$ 19.33  & 82.66 $\pm$ 20.73 \\
BornNet & 61.41 $\pm$ 15.88 & 81.00 $\pm$ 48.87 & 73.74 $\pm$ 18.31 & 78.57 $\pm$ 19.62 \\
LinearBornNet & 50.63 $\pm$ 13.50 & \textbf{59.48 $\pm$ 31.07} & 73.92 $\pm$ 19.11 & 75.60 $\pm$ 18.67 \\
BornAE & \textbf{49.80 $\pm$ 11.78} & 82.44 $\pm$ 29.22 & 76.94 $\pm$ 18.91 & 83.74 $\pm$ 21.63 \\
\hline
Class & \multicolumn{2}{c||}{\textbf{Shepp-Logan}} & \multicolumn{2}{c|}{\textbf{Rectangles}} \\
\hline
Method & $\ell^2$ (\%) & $\ell^1$ (\%) & $\ell^2$ (\%) & $\ell^1$ (\%) \\
\hline
Born ($\gamma=1$) & 97.94 $\pm$ 10.22 & 138.31 $\pm$ 11.24 & 89.53 $\pm$ 24.87 & 118.33 $\pm$ 31.83 \\
NIO               & 75.05 $\pm$  8.61 &  91.70 $\pm$ 11.58 & 57.44 $\pm$  5.08 & 59.36 $\pm$  5.67 \\
BCR-Net           & 78.37 $\pm$  8.03 &  87.31 $\pm$ 10.10 & 73.81 $\pm$  6.44 & 75.09 $\pm$  8.47 \\
CNN               & 59.29 $\pm$  9.68 & 84.09 $\pm$ 12.24 & 42.62 $\pm$  3.84 & 58.05 $\pm$  5.45 \\
CNNB              & 58.59 $\pm$ 8.69 & 87.11 $\pm$ 12.73 & 42.04 $\pm$ 5.89 & 55.17 $\pm$ 8.38 \\
BCNN              & 59.03 $\pm$  8.59 & 102.32 $\pm$ 14.02 & 45.95 $\pm$  9.26 & 75.37 $\pm$ 16.07 \\
BornNet & 57.12 $\pm$ 10.92 & \textbf{58.69 $\pm$ 18.31} & \textbf{33.83 $\pm$ 4.69} & 28.92 $\pm$ 6.73 \\
LinearBornNet & \textbf{52.99 $\pm$ 11.35} & 61.31 $\pm$ 17.61 & 34.92 $\pm$ 6.80 & \textbf{26.89 $\pm$ 13.48}\\
BornAE & 60.67 $\pm$ 9.89 & 84.72 $\pm$ 13.46 & 40.07 $\pm$ 5.94 & 54.21 $\pm$ 11.38 \\
\hline
\end{tabular}
\label{tab:errors}
\end{table}

\begin{figure}[!htbp]
\centering
\includegraphics[width=\linewidth]{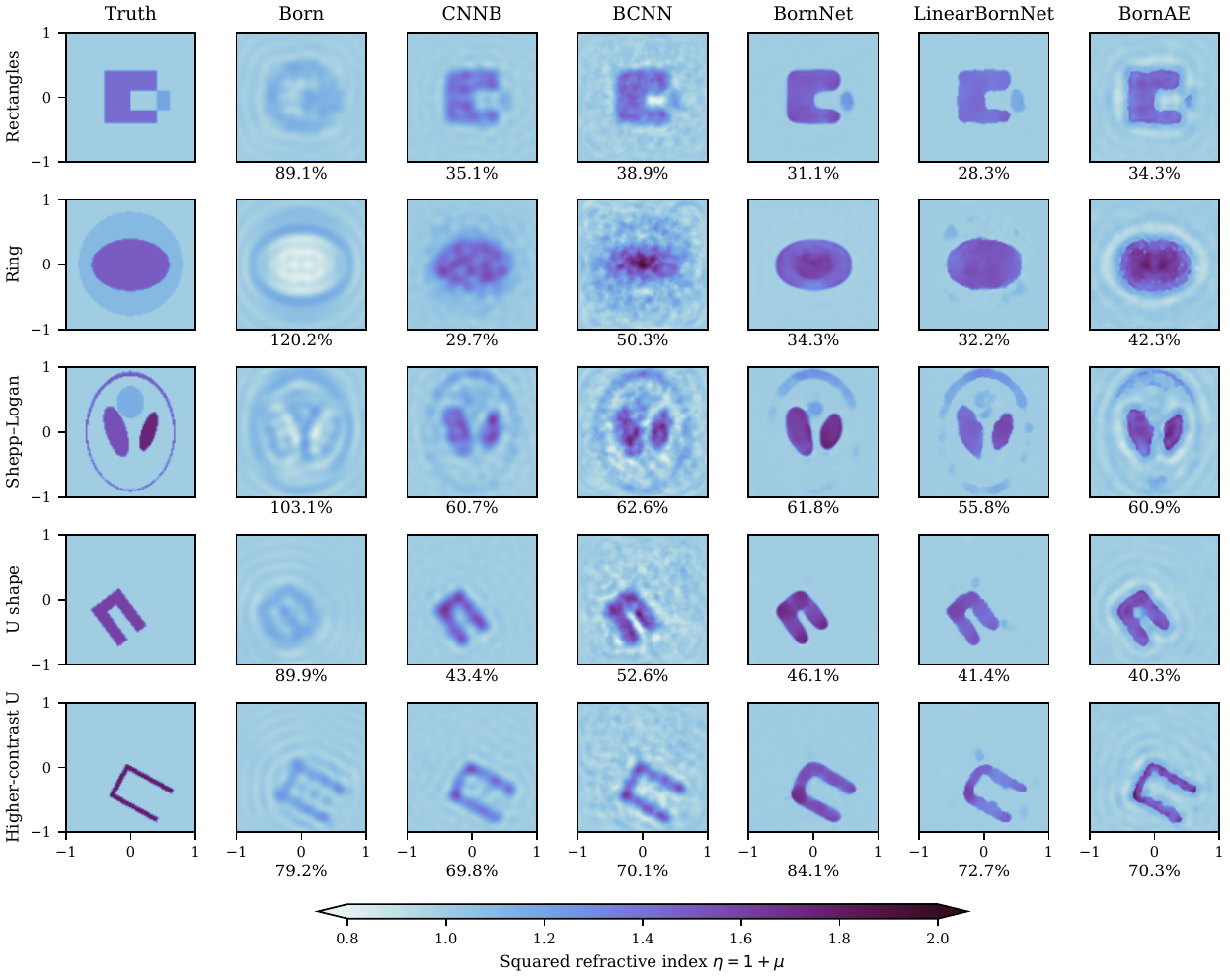}
\caption{Noiseless reconstructions of five penetrable objects with
shapes absent from training.  Columns show the truth and the six
reconstruction methods.  Numbers give $100E_\mu$.  The two U-shaped
objects differ in geometry as well as contrast.}
\label{fig:farfield-ood}
\end{figure}

Table \ref{tab:errors} summarizes error statistics for the four
out-of-distribution classes listed there, and Figure~\ref{fig:farfield-ood} shows the first realization of the test scatterer from
 each of five scatterer configurations.  No noise is added to the data and no network is retrained for these examples.

LinearBornNet has the smallest error on the rectangular example,
while CNNB has the smallest error on the ring.  All methods miss
fine features of the Shepp--Logan phantom, suggesting the need to allow intersection of training disks when we train the networks.  In the neural reconstructions,
the openings of the U-shaped objects remain visible, but their widths
and the wall thicknesses are
inaccurate, particularly for the thin, higher-contrast example.
The relative errors for that example range from 69.8\% to 84.1\%
among the networks.

\subsection{Imaging impenetrable scatterers}
\label{sec:obstacles}
All the neural networks in this paper are trained on data for penetrable (and non-absorbing) media.  In this section we investigate the response of the network to data arising from out of distribution  impenetrable scatterers.  None of the networks were retrained for this example. We compute data for  experiments using one, two or three  impenetrable discs in the measurement area that replace the penetrable inclusions by sound-hard obstacles
occupying a region $D$.  The total field satisfies \eqref{eq:helm_ext}-\eqref{eq:soundhard}.

We generate sound-hard far-field data at $\kappa=16$ for the 4,000 disk
geometries used in the near-field test set in Section~\ref{sec:near}:
1,343 single disks, 1,352 pairs and 1,305 triples. The radii lie in
$[0.1,0.3]$, with centers proposed uniformly in $[-0.7,0.7]^2$ and
overlaps rejected. Here we use plane-wave illumination and the original
far-field acquisition, reconstruction grid and trained networks, without
adding noise or retraining. The fields are computed by coupled
cylindrical-wave expansions, including scattering between the disks.
We increase the truncation order until the relative change in the
far-field matrix is below $10^{-6}$. The retained orders ranges from
40 to 270.

The networks were trained to output a real material contrast, but
this interpretation does not apply to a sound-hard obstacle. Each
network still returns one real-valued image, which may contain
positive and negative values, and it does not predict an imaginary
contrast. We therefore use the output only as an indicator of shape.
For each image, pixels with absolute value at least half its maximum
absolute value define the predicted obstacle mask. We compare this
mask with the true obstacle mask using the intersection-over-union
(IoU) score: the number of pixels in their intersection divided by
the number in their union. IoU is one for identical masks and zero
for disjoint masks.

\begin{table}[!htbp]
\caption{Shape IoU (mean $\pm$ standard deviation) for 4,000 sound-hard
test objects: 1,343 single disks, 1,352 pairs and 1,305 triples.
Larger values indicate better agreement with the true obstacle shapes;
bold entries identify the largest mean IoU in each column.}
\centering
\small
\setlength{\tabcolsep}{4pt}
\begin{tabular}{lrrrr}
\toprule
Method & One disk & Two disks & Three disks & All objects\\
\midrule
Born & $0.345\pm0.067$ & $0.289\pm0.068$ & $0.257\pm0.063$ & $0.297\pm0.076$\\
CNNB & $0.292\pm0.095$ & $0.225\pm0.071$ & $0.194\pm0.062$ & $0.237\pm0.087$\\
BCNN & $0.143\pm0.071$ & $0.132\pm0.050$ & $0.126\pm0.044$ & $0.134\pm0.057$\\
BornNet & $0.468\pm0.093$ & $0.472\pm0.071$ & $0.425\pm0.093$ & $0.455\pm0.089$\\
LinearBornNet & $\mathbf{0.499\pm0.080}$ & $\mathbf{0.475\pm0.056}$ & $\mathbf{0.433\pm0.062}$ & $\mathbf{0.469\pm0.072}$\\
BornAE & $0.494\pm0.068$ & $0.426\pm0.062$ & $0.360\pm0.064$ & $0.427\pm0.085$\\
\bottomrule
\end{tabular}
\label{tab:obstacles}
\end{table}

\begin{figure}[!htbp]
\centering
\includegraphics[width=\linewidth]{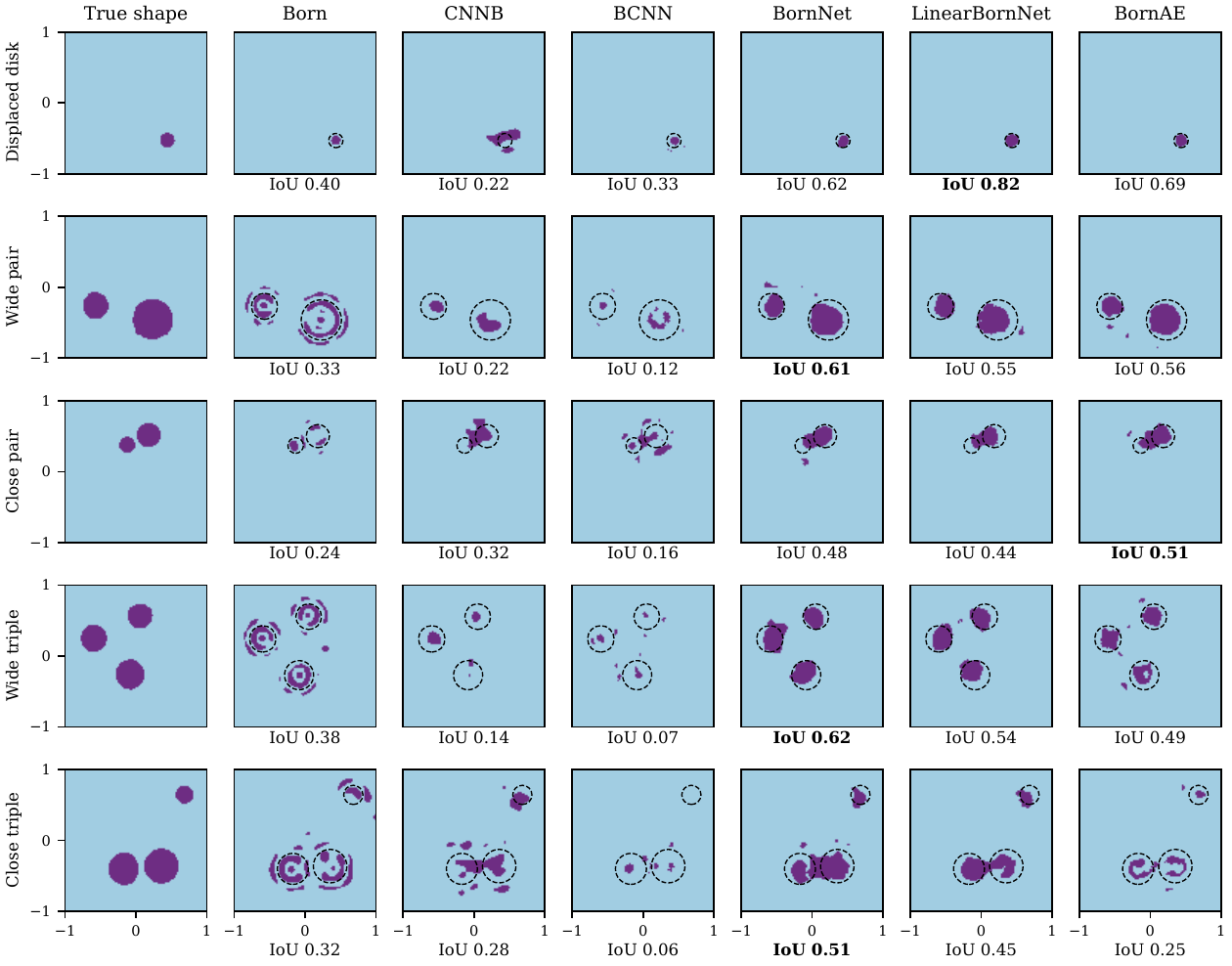}
\caption{Selected sound-hard shape reconstructions from the 4,000-object
noiseless far-field test set.
Dark-purple pixels show the true obstacle mask in the first column and
the predicted masks in the remaining columns. Predictions are
thresholded at half the maximum absolute reconstruction value.
Dashed circles mark the true boundaries. Numbers give shape IoU;
bold values identify the best displayed method for each object.
Each panel covers $\Omega=[-1,1]^2$.}
\label{fig:farfield-obstacles}
\end{figure}

LinearBornNet gives the largest mean IoU overall and for each disk
count (Table~\ref{tab:obstacles}). All three new networks have larger
mean IoUs than the Born reconstruction in each group, while CNNB and
BCNN have smaller means.
Figure~\ref{fig:farfield-obstacles} shows selected reconstructions of a
displaced disk, two pairs and two triples. These examples illustrate
successful localization as well as difficulties with close disks;
the table includes all 4,000 objects.

All three new networks locate the well-separated obstacles, but their
sizes and boundaries remain inaccurate. In the close-pair example,
all three merge the disks. For the close triple, BornNet and
LinearBornNet merge two neighboring disks, while BornAE separates
them but produces fragmented masks. Additional patches also appear
outside the true objects. Nonetheless, the new networks do give positional information for sufficiently well separated impenetrable scatterers, even though they were not part of the training data.

\section{Near-field imaging}\label{sec:near}
In non-destructive testing applications, the acoustic sources and receivers are often near to the scatterer under scrutiny.  The Born approximation is now no longer a band-limited Fourier transform so we are interested in methods like LinearBornNet and BornAE that do not require matrix inversion.

To test the networks on near-field imaging, we use point-source illumination and
measurements of the scattered field at a finite distance from the scatterers.
The source and receiver positions are
$\bfs_j=(R^T\cos\theta_j,R^T\sin\theta_j)$ and
$\bfx_\ell=(R^M\cos\theta_\ell,R^M\sin\theta_\ell)$,
respectively, with $1\leq j,\ell\leq N_m$.
Both circles must be chosen to lie outside the search region $\Omega=[-1,1]^2$.
For the incident field $u^i(\bfy,\bfs_j)=\Phi(\bfy,\bfs_j)$ defined
in Eq.~\eqref{eq:ui_point}, the Born approximation is
\begin{equation}
    u^s(\bfx_\ell,\bfs_j)\approx
\kappa^2\int_\Omega
\Phi(\bfx_\ell,\bfz)\Phi(\bfz,\bfs_j)\mu(\bfz)\,dA_z.
\end{equation}
Unlike the far-field formula, the kernel of this integral is a product of two Green
functions rather than a Fourier kernel.

For the near-field experiments, we use a cell-centered grid with
spacing $h_{\rm NF}=2/N_g$:
\begin{equation}
    \bfz^{\rm NF}_{p,q}
:=\left(-L/2+(p-\tfrac12)h_{\rm NF},
        -L/2+(q-\tfrac12)h_{\rm NF}\right),
\qquad 1\leq p,q\leq N_g.
\end{equation}
Midpoint quadrature defines the discrete near-field Born operator:
\begin{equation}
    (\mathscr{B}_a^{\rm NF}\mu)_{\ell,j}
:=h_{\rm NF}^2\kappa^2\sum_{q,r=1}^{N_g}
\Phi(\bfx_\ell,\bfz^{\rm NF}_{q,r})
\Phi(\bfz^{\rm NF}_{q,r},\bfs_j)\mu(\bfz^{\rm NF}_{q,r}).
\end{equation}
After vectorization,
$\mathscr{B}_a^{\rm NF}\in\mathbb C^{N_m^2\times N_g^2}$. We solve for a real contrast by stacking the real and imaginary parts and performing a Tikhonov regularized inversion as in the far field case."

For the near-field experiments we chose $\gamma=6.54\times10^{-4}$ by finding the value that minimizes the mean relative
contrast error on noise-free validation data. The search for the best regularization parameter uses ten
logarithmically spaced values of $\sqrt\gamma/\|A_{\rm NF}\|_2$
between $10^{-5}$ and $10^{-1/2}$.

The experiments use $\kappa=16$, $L=2$, a $100\times100$ image grid,
and 100 equally spaced source and receiver positions on the same
circle of radius $R^T=R^M=1.8$ (approximately one wavelength from the corners of $\Omega$).

Each training, validation or test scatterer consists of one, two or three non-overlapping disks, with radii
drawn from $\mathcal U(0.1,0.3)$, centers proposed uniformly in
$[-0.7,0.7]^2$, and independent squared refractive indices drawn
from $\mathcal U(1.1,1.8)$ as before.
We generated the scattered fields using coupled
cylindrical-wave expansions with orders $-24,\ldots,24$, rather
than FEM. The measurements contain only the scattered field,
including when the source and receiver coincide.
To represent circular boundaries on the image grid, each pixel is
divided into $8\times8$ equal subcells, and its target value is the
average contrast at their 64 centers. A pixel cut by a disk boundary
therefore approximates an area-weighted average of the contrasts
inside and outside the disk.

The data comprises 16,000 training, 4,000 validation, and 4,000 test
objects. Each training run used 25,000 Adam steps, batch size 16,
learning rate $3\times10^{-4}$, and the same random seed.
We added no measurement noise during training. For each run, model
selection used the lowest mean relative contrast error on the 4,000
noiseless validation objects, evaluated every 1,000 training steps.
Inputs were divided by $s=8.4473\times10^{-3}$, the RMS of the real
and imaginary training channels, and contrast images were not scaled.
Accordingly, BornNet uses $sK_{\rm NF}$ as its fixed inverse layer,
and BornAE uses $\mathscr{B}_a^{\rm NF}/s$ before the learned field
corrector. LinearBornNet uses its learned inverse (i.e. no Born approximation is used).

We also retrained CNNB and BCNN, retaining their previously tuned
CNN architectures. CNNB is trained to predict the scaled Born field
$A_{\rm NF}\mu/s$; applying $sK_{\rm NF}$ gives the contrast
reconstruction. BCNN learns an additive correction to the Born
reconstruction $K_{\rm NF}x_{\rm NF}$, minimizing mean squared error
between the final reconstruction and the true contrast.

For the reported near-field BornAE model, we use the same three loss
terms as in Eq.~\eqref{eq:bornae-loss}, evaluating the field error
after undoing the normalization. Let $D_\phi^{\rm NF}(\widehat\mu)$
denote the decoder output in normalized field units. The corresponding
un-normalized prediction is $sD_\phi^{\rm NF}(\widehat\mu)$, so the loss is
\begin{equation}
 \mathcal L_{\rm NF}:=
 \frac{\|\widehat\mu-\mu\|_2^2}{N_g^2}
 +\frac{1}{2N^MN^T}
   \|sD_\phi^{\rm NF}(\widehat\mu)-x_{\rm NF}\|_2^2
 +\frac{10^{-4}}{N_g^2}\|\widehat\mu\|_1.
 \label{eq:bornae-loss-nf}
\end{equation}
Equivalently, the field term is $s^2$ times the MSE between the
normalized decoder output and $x_{\rm NF}/s$. This compensates for
the factor $1/s^2$ introduced by normalizing the fields, retaining
unit weight on the un-normalized field MSE as in the far-field loss.
The loss structure is therefore the same, while the near-field
training settings and model-selection procedure are those given above.

An initial run instead used the normalized field MSE without this
compensation and gave a noiseless mean test error of 119.08\%.
We then compared two separate training runs multiplying that
normalized field MSE by $10^{-2}$ and $s^2$, respectively, with
the other settings fixed. Clean validation error selected the latter choice of $s$,
corresponding to Eq.~\eqref{eq:bornae-loss-nf}; all near-field BornAE
results below use this model.

While testing, each complex measurement is perturbed
as in Eq.~\eqref{eq:measurement-noise}.
For each test object, one noise array is reused for all methods and
all noise levels; changing $\delta$ changes only its amplitude.
Table~\ref{tab:nearfield-minimal} reports the contrast error from
Eq.~\eqref{eq:contrast-error}.

\begin{table}[t]
\caption{Near-field relative contrast errors $100E_\mu$: mean $\pm$
standard deviation across 4,000 test objects. Bold entries mark the
lowest mean error at each noise level. BornAE uses the loss
in Eq.~\eqref{eq:bornae-loss-nf}.}
\centering
\small
\setlength{\tabcolsep}{3pt}
\begin{tabular}{lrrrrr}
\toprule
Method & $\delta=0$ & $\delta=0.1$ & $\delta=0.25$ & $\delta=0.5$ & $\delta=1$\\
\midrule
Born & $66.12\pm20.32$ & $66.14\pm20.31$ & $66.23\pm20.27$ & $66.58\pm20.13$ & $67.93\pm19.61$\\
CNNB & $31.34\pm5.20$ & $31.43\pm5.20$ & $33.35\pm5.33$ & $42.55\pm6.29$ & $66.53\pm8.44$\\
BCNN & $34.42\pm5.14$ & $34.98\pm5.03$ & $39.11\pm5.06$ & $52.60\pm10.15$ & $81.33\pm20.35$\\
BornNet & $\mathbf{5.71\pm2.89}$ & $\mathbf{5.82\pm2.86}$ & $\mathbf{6.54\pm2.88}$ & $\mathbf{9.22\pm3.26}$ & $\mathbf{18.83\pm5.93}$\\
LinearBornNet & $8.12\pm4.03$ & $11.63\pm7.19$ & $25.50\pm21.26$ & $53.47\pm42.97$ & $94.94\pm76.02$\\
BornAE & $10.51\pm3.83$ & $10.56\pm3.83$ & $10.97\pm3.81$ & $12.75\pm3.85$ & $20.57\pm5.00$\\
\bottomrule
\end{tabular}
\label{tab:nearfield-minimal}
\end{table}

From Table~\ref{tab:nearfield-minimal}, BornNet has the lowest mean error at every noise level.
CNNB and BCNN give clean-data errors of 31.34\% and 34.42\%,
rising to 66.53\% and 81.33\% at $\delta=1$, respectively.
LinearBornNet achieves an 8.12\% error on clean data, but this rises
to 94.94\% at $\delta=1$.
Using the loss in Eq.~\eqref{eq:bornae-loss-nf}, BornAE gives errors of
10.51\% on clean data and 20.57\% at $\delta=1$, compared with
5.71\% and 18.83\% for BornNet.

\begin{figure}[!htbp]
\centering
\includegraphics[width=\linewidth]{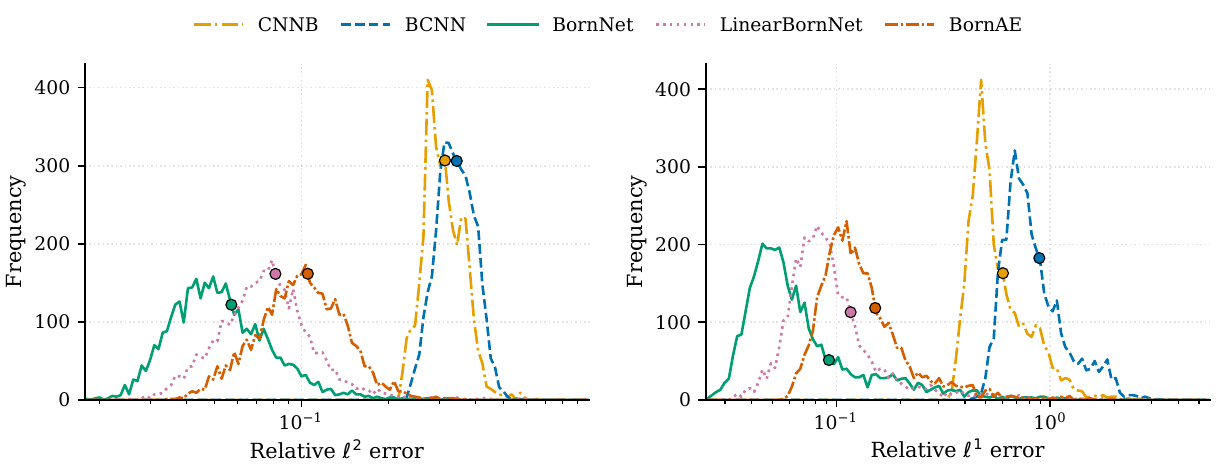}
\caption{Relative contrast errors in the $\ell^2$ ($E_{\mu,2}$) (left) and $\ell^1$ ($E_{\mu,1})$
(right) norms for 4,000 noiseless near-field test objects. Curves show
histogram counts with common logarithmically spaced bins within each
panel; circles mark the means, as in Fig.~\ref{fig:farfield-errors}.
Logarithmic horizontal axes display the full error range.
BornAE uses the same model as in Table~\ref{tab:nearfield-minimal}.}
\label{fig:nearfield-errors}
\end{figure}

Figure~\ref{fig:nearfield-errors} compares the clean-data ($\delta=0$) error
distributions.  Compared with the far-field results in Fig.~\ref{fig:farfield-errors}, the near-field results for BornNet and LinearBornNet no longer give essentially equivalent results.  BornNet has an advantage both in average error and in the error distribution.  

\begin{figure}[!htbp]
\centering
\includegraphics[width=\linewidth]{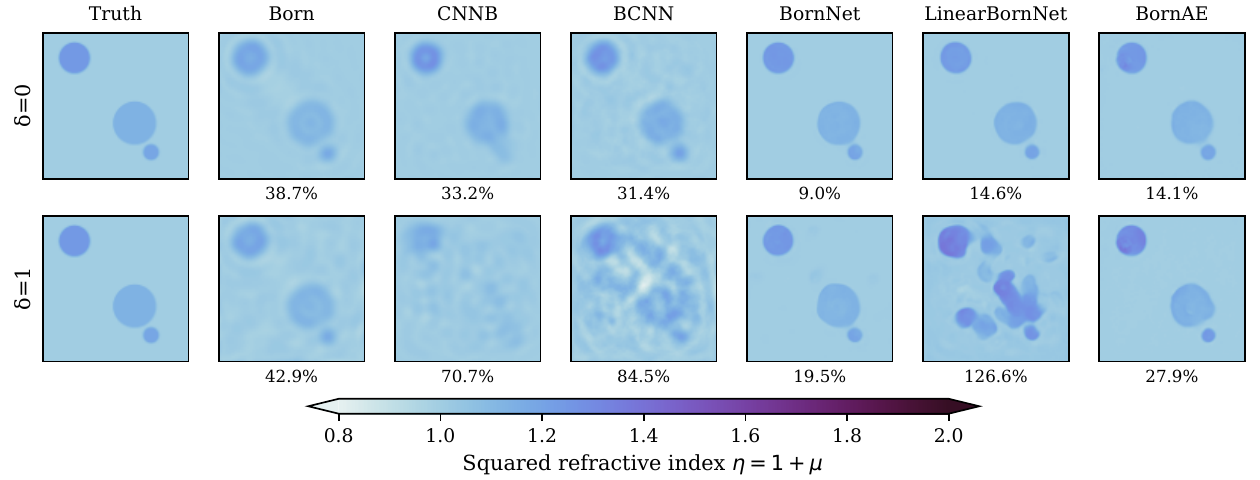}
\caption{Near-field reconstructions of the first three-disk test object
in the test data. The top row uses clean data; the bottom row uses
the same noisy measurements ($\delta=1$) for every method.
BornAE uses the same model as in Table~\ref{tab:nearfield-minimal},
trained with the loss in Eq.~\eqref{eq:bornae-loss-nf}.
Colors show $n^2=1+\mu$ using the same color map and scale $[0.8,2]$
as Fig.~\ref{fig:farfield-ood};
values outside this range take the endpoint colors.
Numbers give $100E_\mu$, computed from the unclipped reconstructions.}
\label{fig:nearfield-minimal}
\end{figure}

In Fig.~\ref{fig:nearfield-minimal} we show a realization of a test example consisting of three scatterers. All methods (including Born - this is a weak scattering example) reconstruct the scatterers in the noise-free case.  Surprisingly, CNNB performs poorly at higher noise levels.  Both BornNet and BornAE perform well at higher noise levels where both reconstruct the three disks.
LinearBornNet introduces large hallucinations outside the true disks.

\section{Conclusions}\label{sec:conc}

We have proposed, trained, and tested three new neural networks for inverse scattering.
They are motivated by the goal of reducing error compared to our previous neural
networks, and investigating if the Born approximation is even necessary
to improve the network performance.

BornNet and LinearBornNet have similar mean errors on noiseless
far-field data, while LinearBornNet has lower mean errors at every
tested far-field noise level. In near-field tests, BornNet has
the lowest mean error at all noise levels.

In 3D it could be more efficient to avoid the need to invert the regularized Born operator due to the size of the corresponding full matrix for the Born method. In this regard BornAE is promising, since it only uses the forward Born matrix during training. Its mean error increases less from noise-free data ($\delta=0$) to $\delta=1$ than that of the other two new methods in both settings.

Necessarily, a study such as this must be speculative, since the results
are only computed for a limited selection of targets.  This is particularly true for the near-field study that only has examined in-distribution targets (albeit with noise).  Obviously what is needed is a theoretical basis for choosing between the neural networks.  

Future directions include the use of  LinearBornNet and BornAE in three dimensional
inverse scattering, particularly when trained on data including Shepp-Logan type
scatterers.  Another interesting study would be to apply the scheme to measured data.

\section*{Acknowledgments}

The research of T.L. is supported by the Research Council of Finland via the research project 373639 and the Flagship of Advanced Mathematics for Sensing Imaging and Modeling grant 358944. We also thank Professor
Ke Chen of the University of Delaware for sharing his computational resources.

\section*{AI Acknowledgment}
The authors used generative-AI tools (OpenAI’s ChatGPT) during the development and preparation of this paper to assist with code development and computational experiments, improve exposition, and check for typographical and grammatical errors. This assistance supplemented the authors’ own coding, experimentation, and editing. AI-generated outputs were treated as unverified suggestions and reviewed by the authors. The authors made all final decisions concerning the code, experiments, and manuscript and take full responsibility for the paper’s contents.

\end{document}